\documentclass[11pt,letterpaper,reqno]{amsart}
\usepackage[
 top=0.9in,
 bottom=0.65in,
 inner=0.95in,
 outer=0.95in,
 marginparwidth=2cm
]{geometry}
\usepackage[T1]{fontenc}
\usepackage{microtype}
\usepackage{amsmath,amssymb,amsthm,mathtools,booktabs}
\usepackage{placeins}
\usepackage{multirow}
\usepackage{array}
\usepackage{longtable}
\usepackage{listings}
\usepackage{tikz}
\usetikzlibrary{arrows.meta,calc}
\usepackage{hyperref}
\hypersetup{colorlinks=false,
 pdfborder={0 0 1},
 linkbordercolor={1 0 0},
 citebordercolor={0 0 1},
 urlbordercolor={0 0 1},
 filebordercolor={0 0 1},
 bookmarksnumbered=true,
 bookmarksdepth=2,
 pdftitle={Higher-Derivative Conjectures of Sendov},
 pdfsubject={Analytic proofs of Conjectures 3, 8, and 9},
 pdfkeywords={Sendov conjecture, higher derivatives, polynomial zeros, quartic polynomials}}

\allowdisplaybreaks[4]
\makeatletter
\def\th@plain{\thm@notefont{}\itshape}
\def\th@definition{\thm@notefont{}\normalfont}
\def\subsubsection{\@startsection{subsubsection}{3}%
  \z@{.5\linespacing\@plus.7\linespacing}{-.5em}%
  {\normalfont\bfseries\let\@secnumfont\bfseries}}
\makeatother
\numberwithin{equation}{section}
\numberwithin{table}{section}
\numberwithin{figure}{section}
\newtheorem{mainthm}{Theorem}

\newtheorem{theorem}{Theorem}[section]
\newtheorem{lemma}[theorem]{Lemma}
\newtheorem{proposition}[theorem]{Proposition}
\newtheorem{corollary}[theorem]{Corollary}
\theoremstyle{remark}\newtheorem{remark}[theorem]{Remark}
\newtheorem*{remark*}{Remark}
\newcommand{\D}{\overline{\mathbb D}}
\newcommand{\E}{\mathbb E}
\newcommand{\N}{\mathbb N}
\newcommand{\Rea}{\operatorname{Re}}
\newcommand{\Sfun}{S_s}
\newcommand{\C}{\mathbb C}
\newcommand{\R}{\mathbb R}
\newcommand{\I}{\mathbf{i}}
\newcommand{\Bbern}{\mathcal B}
\providecommand{\abs}[1]{\lvert#1\rvert}
\providecommand{\Absbig}[1]{\bigl\lvert#1\bigr\rvert}
\providecommand{\Absbigg}[1]{\biggl\lvert#1\biggr\rvert}

\providecommand{\AbsBigg}[1]{\Biggl\lvert#1\Biggr\rvert}

\renewcommand{\=}{\coloneqq}
\renewcommand{\le}{\leqslant}
\renewcommand{\leq}{\leqslant}
\renewcommand{\ge}{\geqslant}

\newenvironment{appendixtable}{%
 \begingroup\small
 \renewcommand{\arraystretch}{1.3}%
 \setlength{\extrarowheight}{2pt}%
 \setlength{\tabcolsep}{7pt}%
 \setlength{\heavyrulewidth}{0.6pt}%
 \setlength{\lightrulewidth}{0.3pt}%
 \setlength{\LTleft}{0pt}%
 \setlength{\LTright}{0pt}%
 \setlength{\LTcapwidth}{\linewidth}%
 \setlength{\LTpre}{5pt plus 2pt minus 1pt}%
 \setlength{\LTpost}{4pt plus 2pt minus 1pt}%
}{\endgroup}
\title[Higher-Derivative Conjectures of Sendov]{Higher-Derivative Conjectures of Sendov}
\makeatletter
\renewcommand{\@setauthors}{\begingroup
	\trivlist
	\centering\normalsize \@topsep30\p@\relax
	\advance\@topsep by -\baselineskip
	\item\relax
	\authors
	\vspace{14pt}
	\endtrivlist
	\endgroup
}
\makeatother
\author[Mukhiddin~Ismoilov, Zhiqiang~Li, Karim~Rakhimov, and Chengliang~Tan]{
\resizebox{\textwidth}{!}{
    \begin{tabular}[t]{c}Mukhiddin~Ismoilov\\[-3pt]\Small International Islamic \\ [-3pt]\small Academy of Uzbekistan\\{ } \end{tabular}\;
	\begin{tabular}[t]{c}Zhiqiang~Li\\[-3pt]\Small Peking University\\{ }  \end{tabular}\;
	\begin{tabular}[t]{c}Karim~Rakhimov\\[-3pt]\Small V.~I.~Romanovsky Institute of Mathematics,\\ [-3pt]\small Uzbekistan Academy of Sciences \end{tabular}\;
	\begin{tabular}[t]{c}Chengliang~Tan\\[-3pt]\Small Peking University \\ { } \end{tabular}
}
}

\begin{document}
\raggedbottom
\begin{abstract}
We prove Conjectures~3 and~8 of Sendov's 2002 survey and deduce Conjecture~9. More precisely, let $p$ be a polynomial of degree $n$ whose zeros lie in the closed unit disk. For every integer $s$ with $1\le s<n$, we prove that each zero of $p$ lies within distance $2s/(s+1)$ of a zero of $p^{(s)}$. For quartics, we obtain the sharp bound $2/\sqrt3$ for the second derivative and determine all polynomials attaining it.
\end{abstract}
\maketitle

\section{Introduction and main results}\label{sec:mainresults}

The relation between the zeros of a polynomial and the zeros of its derivatives is a central topic in the geometry of polynomials~\cite{Marden1966,RahmanSchmeisser2002,SheilSmall2002}. In particular, the critical points are precisely the points at which a nonconstant polynomial fails to be locally injective, so their location is relevant to the geometry of polynomials. For polynomials with real coefficients, Rolle's theorem asserts that there is at least one real critical point between any two consecutive distinct real zeros. For complex polynomials, the Gauss--Lucas theorem says that all critical points are in the convex hull of the zeros. These geometric conclusions suggest a more precise quantitative question: how close must a critical point be to each prescribed zero? Grace and Heawood showed that if a complex polynomial of degree $n\ge2$ has two distinct zeros $a$ and $b$, then its derivative has a zero in the closed disk with center $(a+b)/2$ and radius $\frac12\abs{a-b}\cot(\pi/n)$~\cite{Marden1985}. Sendov made significant contributions to this area by establishing distance bounds in several cases and formulating a number of conjectures, including the well-known Sendov conjecture which was proposed in 1958. He also presented some of these results and conjectures in a survey~\cite{Sendov2002}.

To state these conjectures, we first introduce some notation. Set $\N\=\{1,\,2,\,\ldots\}$ and $\N_0\=\{0\}\cup \N$. For a nonconstant polynomial $p\in\C[z]$, write $A(p)$ for its set of distinct zeros and $H(p)\=\operatorname{conv}A(p)$ for the convex hull of $A(p)$. For nonempty compact subsets $E,\,F$ of $\C$, define
\begin{equation*}
 \rho(E,F)\=\max_{x\in E}\min_{y\in F}\abs{x-y}.
\end{equation*}
Thus, $\rho(E,F)$ is the maximum distance from a point of $E$ to the set $F$ and is the main quantity studied in this article. It is worth noting that $\rho(E,F)\ne\rho(F,E)$ in general. Let $\D\=\{z\in\C:\abs{z}\le1\}$ and $p^{(s)}$ be the $s$-th derivative of $p$. If $\deg p=n\ge2$ and $A(p)\subseteq\D$, the Gauss--Lucas theorem gives $A\bigl(p^{(s)}\bigr)\subseteq H(p)$ and hence $\rho\bigl(A(p),A\bigl(p^{(s)}\bigr)\bigr)\le2$ for $1\le s<n$. The classical Sendov conjecture asserts that $\rho(A(p),A(p'))\le1$, where the upper bound is optimal for every degree, as shown by $p(z)=z^n-1$.

Sendov's conjecture remained open for nearly seven decades and was recently resolved by an AI-assisted proof~\cite{Mazur2026,Tao2026}. Before the conjecture was finally resolved, many mathematicians had obtained important partial results. Brown and Xiang verified Sendov's conjecture for polynomials with at most eight distinct zeros using upper and lower bounds for products of the moduli of critical points~\cite{BrownXiang1999}. D\'egot showed that the conjecture holds at a prescribed zero $a\in(0,1)$ for all sufficiently large degrees, with a threshold depending on $a$~\cite{Degot2014}. Chalebgwa made D\'egot's estimates explicit and obtained a uniform threshold when $a$ ranges over a compact subinterval of $(0,1)$~\cite{Chalebgwa2020}. Tao proved the conjecture for all sufficiently large degrees with an absolute threshold independent of $a\in[0,1]$. His proof relies in part on compactness methods and the argument principle~\cite{Tao2022}. Finally, Mazur obtained an AI-assisted proof in all degrees, based on integral identities and exact rational verification, together with a separate Lean formalization of the theorem~\cite{Mazur2026}. Tao subsequently gave a simplified proof based on polynomial integral identities and Maclaurin's inequality, together with a Lean formalization developed with AI assistance. In fact, the proof gives the stronger conclusion conjectured by Phelps and Rodriguez: if $p$ is a complex polynomial of degree $n\ge2$ whose zeros lie in $\D$, then for every zero $a$ of $p$, there exists a zero $\zeta$ of $p'$ satisfying $\abs{\zeta-a}<1$, unless $\abs{a}=1$ and $p(z)=c(z^n-a^n)$ for some nonzero $c\in\C$~\cite{Tao2026}.

Beyond the case of critical points, it is also natural to study the distances between the zeros of a polynomial and the zeros of its higher derivatives. Let $p\in\C[z]$ satisfy $\deg p=n\ge2$ and $A(p)\subseteq\D$, and let $s$ be an integer with $2\le s<n$. Repeated application of the classical Sendov conjecture, now a theorem, gives $\rho\bigl(A(p),A\bigl(p^{(s)}\bigr)\bigr)\le s$. This estimate does not even yield a bound smaller than $2$, which follows directly from the Gauss--Lucas theorem. Sendov also considered this problem and proposed the following conjecture, which is formulated as our first main result in this article.

\begin{mainthm}[{\cite[Conjecture~3]{Sendov2002}}]\label{thm:main}
Let $p\in\C[z]$ satisfy $\deg p=n\ge2$ and $A(p)\subseteq\D$.
Then, for every integer $1\le s<n$,
\begin{equation*}\rho\bigl(A(p),A\bigl(p^{(s)}\bigr)\bigr)\le\frac{2s}{s+1}.\end{equation*}
\end{mainthm}

\begin{remark}\label{rem:generalsharp}
Sendov stated the conjecture for monic polynomials $p(z)$, for which $\D$ is the smallest disk containing all zeros of $p$, and the present formulation is equivalent by the usual normalization. 

For each fixed $s\in\N$, the constant $2s/(s+1)$ in the theorem is optimal uniformly over all degrees $n>s$. Indeed, $p_s(z)\=(z-1)(z+1)^s$ has degree $s+1$ and all its zeros are in $\D$. Since $p_s^{(s)}(z)=s!((s+1)z+s-1)$, we have $\rho\bigl(A(p_s),A\bigl(p_s^{(s)}\bigr)\bigr)=1+\frac{s-1}{s+1}=\frac{2s}{s+1}$. 

For a general polynomial $p$ satisfying the hypotheses of Theorem~\ref{thm:main} with $n=s+1$, a simple calculation shows that $\rho\bigl(A(p),A\bigl(p^{(s)}\bigr)\bigr)=\frac{2s}{s+1}$ if and only if $p(z)=c(z-a)(z+a)^{n-1}$ for some nonzero $c\in\C$ and $a\in\C$ with $\abs{a}=1$.
\end{remark}

Sendov obtained the upper bound at a prescribed zero on the unit circle, for $n-s=2,\,3$, and for $n-s=4$ when $n\ge6$. He also proved it when all vertices of the convex hull of the zeros lie on the boundary of the smallest closed disk containing all the zeros~\cite{Sendov2002}. Uluchev obtained numerical evidence which supports the upper bound $2s/(s+1)$ for several pairs $(n,s)$~\cite{Uluchev2003}.

In our proof, we first note that, for each zero $a$ of $p$, we can write $p(z)\=(z-a)q(z)$ with $q\in\C[z]$. For $2\le s<n$, set $g\=q^{(s-1)}$. Then
\begin{equation*}
 p^{(s)}(z)=(z-a)g'(z)+sg(z).
\end{equation*}
Here $\deg g=n-s\ge1$, and the Gauss--Lucas theorem gives $A(g)\subseteq\D$. Thus, it suffices to find a zero of $(z-a)g'(z)+sg(z)$ within distance $2s/(s+1)$ of $a$. More precisely, we deduce the case $s\ge2$ of Theorem~\ref{thm:main} from the following theorem. The case $s=1$ follows from the classical Sendov theorem.

\begin{theorem}\label{thm:completeweighted}
Let $s,\,r\in\N$ and $s\ge2$, and let $g\in\C[z]$ satisfy $\deg g=r$ and $A(g)\subseteq\D$. Then for each $a\in\D$,
$G(z)\=(z-a)g'(z)+sg(z)$ has a zero $\zeta$ satisfying
\begin{equation*}
\abs{a-\zeta}\le\frac{2s}{s+1}.
\end{equation*}
\end{theorem}

Our strategy for proving Theorem~\ref{thm:completeweighted} is to argue by contradiction and derive incompatible inequalities from integral identities. More precisely, after a rotation and normalization, we may assume that $a\in[0,1]$ and $g$ is monic. Suppose, for a contradiction, that every zero $\zeta$ of $G$ satisfies $\abs{a-\zeta}>2s/(s+1)$.

We first derive a necessary condition on $a$ for the existence of a counterexample (Lemma~\ref{lem:reduction}). With this restriction on $a$, Lemma~\ref{lem:identities} provides weighted integral representations along line segments starting at $a$. In particular, the two identities in~\eqref{eq:values} are weighted versions of Tao's first origin and polar identities, with weight $t^{s-1}$~\cite[Lemma~6]{Tao2026}. The first inequality in~\eqref{eq:integralbounds} is then used to show that any counterexample must satisfy $U=M+E\ge1$ (Lemma~\ref{lem:originineq}), where $M$ bounds the terms involving $g(0)$ and $G(0)$, while $E$ bounds the integral remainder along the segment from $a$ to $0$ (see~\eqref{eq:U}). To estimate these quantities, we derive upper bounds for the quadratic polynomial $\beta$ defined in~\eqref{eq:betadefinition}. At $s=1$, the formula for $\beta$ reduces to $1-2axt+a^2t^2$, which is exactly the quadratic polynomial in Tao's polar inequality~\cite[(14) and Proposition~10]{Tao2026}.

The estimates for $M$ and $E$ depend on the values of $s$ and $r$. We therefore divide the parameter pairs $(s,r)$ into the four ranges listed in Table~\ref{tab:cases}. In Case~I, the lower bound of $1$ in the second inequality of~\eqref{eq:integralbounds} is incompatible with an upper bound strictly below $1$ for the same weighted integral. For Cases~II--IV, we prove $U<1$ by estimating $M$ and $E$ separately, contradicting the inequality $U\ge1$ in Lemma~\ref{lem:originineq}. The first step in all three cases is to derive an upper bound for $\beta(t)$ before estimating the integral defining $E$. In Case~II, where $s$ is sufficiently large, the bound for $\beta(t)$ can be chosen independently of $s$. Case~III treats $2\le s\le8$ and sufficiently large $r$. Here the estimate for $\beta(t)$ depends on $s$, and the corresponding upper bound for $E$ decreases as $r$ increases. Finally, Case~IV contains finitely many parameter pairs and requires sharper estimates for $\beta(t)$. These arguments exclude every counterexample and complete the proof.

\smallskip

Since $A\bigl(p^{(s)}\bigr)\subseteq H\bigl(p^{(s)}\bigr)$, Theorem~\ref{thm:main} immediately yields the following corollary, which was also conjectured by Sendov.

\begin{corollary}[{\cite[Conjecture~9]{Sendov2002}}]\label{cor:hull}
Let $p\in\C[z]$ satisfy $\deg p=n\ge2$ and $A(p)\subseteq\D$.
Then, for every integer $1\le s<n$,
\begin{equation*}
 \rho\bigl(A(p),H\bigl(p^{(s)}\bigr)\bigr)\le\frac{2s}{s+1}.
\end{equation*}
\end{corollary}

\smallskip

For a fixed degree and derivative order, the optimal constant can be smaller than the uniform bound in Theorem~\ref{thm:main}. For quartic polynomials, Sendov considered the example $f_0(z)\=(z-1)^2\bigl(z^2+\frac23z+1\bigr)$, whose zeros lie on the unit circle. It satisfies $\rho(A(f_0),A(f_0''))=2/\sqrt3$. Motivated by this example and his proof of the bound for quartics with real coefficients~\cite[Theorem~10]{Sendov2002}, Sendov proposed the following refinement. We prove this conjecture in our second main theorem and give a complete classification of equality cases.

\begin{mainthm}[{\cite[Conjecture~8]{Sendov2002}}]\label{thm:quartic}
Let $f\in\C[z]$ satisfy $\deg f=4$ and $A(f)\subseteq\D$.
Then
\begin{equation*}
 \rho(A(f),A(f''))\le\frac2{\sqrt3}.
\end{equation*}
Equality holds if and only if, for some $c\in\C\smallsetminus\{0\}$ and $\theta\in\R$,
\begin{equation}\label{eq:quartic-extremal}
 f(z)=c\bigl(z-e^{\I\theta}\bigr)^2
 \biggl(z^2+\frac23e^{\I\theta}z+e^{2\I\theta}\biggr).
\end{equation}
\end{mainthm}

To prove Theorem~\ref{thm:quartic}, we first introduce some notation. For $n\in \N$ and a set $Z=\{z_1,\,z_2,\,\ldots,\,z_n\}\subseteq\C$, let $\mathcal B(Z)$ denote the smallest closed disk containing $Z$. We also write $\mathcal B(Z)\=\mathcal B(z_1,\,z_2,\,\ldots,\,z_n)$ for brevity. The \emph{radius} of $Z$ is given by
\begin{equation*}
 \operatorname{rad}Z\=\min_{c\in\C}\max_{z\in Z}\abs{z-c}.
\end{equation*}
The minimum exists. Indeed, the maximum is continuous in $c$ and tends to infinity as $\abs{c}\to+\infty$. We then apply the following auxiliary theorem to prove Theorem~\ref{thm:quartic}.

\begin{theorem}\label{prop:fourpoints}
Let $\lambda,\,w_2,\,w_3,\,w_4\in\C$ and $\abs{\lambda}\ge1$ satisfy
\begin{equation}\label{eq:quartic-sums}
 w_2+w_3+w_4=2(1+\lambda)\quad\text{ and }\quad
 w_2w_3+w_2w_4+w_3w_4=6\lambda.
\end{equation}
Then $\operatorname{rad}\{0,\,w_2,\,w_3,\,w_4\}\ge\sqrt3/2$. Equality holds if and only if $\lambda=1$ and $w_2w_3w_4=(88\pm4\sqrt2\,\I)/27$.
\end{theorem}

Our strategy for proving Theorem~\ref{prop:fourpoints} is to reduce the radius estimate to algebraic inequalities in $\lambda$ through an ellipse representation of the three points $w_2,\,w_3,\,w_4$. Write $R\=\operatorname{rad}\{0,\,w_2,\,w_3,\,w_4\}\,$. It suffices to assume $R^2\le3/4$ and show that $R^2=3/4$.

After a translation and a rotation of the four points $0,\,w_2,\,w_3,\,w_4$, the images $\tau_2,\,\tau_3,\,\tau_4$ of $w_2,\,w_3,\,w_4$ lie on a possibly degenerate ellipse centered at the origin. These transformations preserve $R$. We then obtain three lower bounds for $R^2$ (Lemma~\ref{lem:quartic-ellipsebounds}). Algebraic inequalities in $\lambda$ allow us to simplify these estimates and combine them with positive weights (Lemma~\ref{lem:quartic-algebra}). This weighted average is bounded below by $3/4$ for every choice of the disk center, which proves $R^2\ge3/4$ (Lemma~\ref{lem:quartic-elimination}). If $R^2=3/4$, then $\lambda=1$ and the ellipse degenerates to a line segment. The equality conditions then determine the extremal polynomials in Theorem~\ref{thm:quartic}.

The article is organized as follows. Section~\ref{sec:reductions} presents the reductions to Theorems~\ref{thm:completeweighted} and~\ref{prop:fourpoints}, which are proved in Sections~\ref{sec:generalproof} and~\ref{sec:quartic}, respectively. Appendices~\ref{app:finiteidentities}--\ref{app:quartic} record the rational estimates and polynomial identities used in the proofs, and Appendix~\ref{app:code} provides code for reproducing these calculations.

\medskip
\noindent{\bf Acknowledgments.} 
Z.~L.\ and C.~T.\ were partially supported by Beijing Natural Science Foundation (JQ25001) and National Natural Science Foundation of China (12471083). The authors gratefully acknowledge Urgench State University and Peking University for the hospitality during the First and Second Sino--Uzbek Conferences on Complex Analysis in 2025 and 2026. This article grew out of collaborations initiated and developed during these events.

\medskip
\noindent{\bf Generative AI statement.}
ChatGPT was used to assist with drafting some of the proofs and improving the language of the manuscript, which was carefully checked and revised by the authors. It was also used to generate the Python code for the exact computations in this article. The code was thoroughly checked and revised, and the computational results were verified. The authors take full responsibility for the content of this article.

\section{Reduction to two theorems}\label{sec:reductions}
In this section, we give the detailed deductions of Theorems~\ref{thm:main} and~\ref{thm:quartic} from Theorems~\ref{thm:completeweighted} and~\ref{prop:fourpoints}.
\subsection{Reduction of Theorem~\ref{thm:main}}

\begin{proof}[\textbf{Proof of Theorem~\ref{thm:main} from
Theorem~\ref{thm:completeweighted}}]
The case $s=1$ is the classical Sendov theorem~\cite{Mazur2026,Tao2026}. We may therefore assume that $s\ge2$. Fix an arbitrary $a\in A(p)$, write $p(z)=(z-a)q(z)$, and set $g\=q^{(s-1)}$ and $r\=n-s\ge 1$. Then $\deg g=r$ and $A(g)\subseteq\D$ by the Gauss--Lucas theorem.
By Theorem~\ref{thm:completeweighted}, there exists a zero $\zeta$ of
\begin{equation*}
 p^{(s)}(z)=(z-a)q^{(s)}(z)+s q^{(s-1)}(z)=(z-a)g'(z)+sg(z)
\end{equation*}
satisfying $\abs{a-\zeta}\le 2s/(s+1)$.
Since $a\in A(p)$ is arbitrary, the result follows.
\end{proof}
\smallskip

\subsection{Reduction of Theorem~\ref{thm:quartic}}\label{sec:quartic-reduction}

\begin{proof}[\textbf{Proof of Theorem~\ref{thm:quartic} from
Theorem~\ref{prop:fourpoints}}]
Fix $a\in A(f)$.
By a rotation and normalization, we may assume that $a\in[0,1]$ and $f$ is monic.

Let $\zeta_1$ be a zero of $f''$ nearest to $a$.
If $\zeta_1=a$, the desired bound is immediate.
Otherwise, write
\begin{equation*}
 F(w)\=\frac{f(a+(\zeta_1-a)w)}{(\zeta_1-a)^4}
       =:w^4-e_1w^3+e_2w^2-e_3w.
\end{equation*}
Since $F''(1)=f''(\zeta_1)\big/(\zeta_1-a)^2=0$, let $\lambda$ be the other zero of $F''$.
Then $\abs{\lambda}\ge1$ by the choice of $\zeta_1$, and $e_1=2(1+\lambda)$ and $e_2=6\lambda$.
Write the zeros of $F$ as $0,\,w_2,\,w_3,\,w_4$, counted with multiplicity.
By Theorem~\ref{prop:fourpoints}, $\sqrt3/2\leq\operatorname{rad}A(F)=\operatorname{rad}A(f)/\abs{\zeta_1-a} \le1/\abs{\zeta_1-a}$, which implies that $\abs{\zeta_1-a}\le2/\sqrt3$ for every $a\in A(f)$.

For $f$ given by~\eqref{eq:quartic-extremal}, direct calculation gives $f''(z)=12c\bigl(z-(1/3)e^{\I \theta}\bigr)^2$ and $\rho(A(f),A(f''))=2/\sqrt3$.
Conversely, suppose equality holds in Theorem~\ref{thm:quartic} and choose $a\in A(f)$ with $\abs{\zeta_1-a}=2/\sqrt3$.
Then $\sqrt3/2\le\operatorname{rad}A(F) \le1/\abs{\zeta_1-a}=\sqrt3/2$, so $\operatorname{rad}A(F)=\sqrt3/2$.
Theorem~\ref{prop:fourpoints} gives $\lambda=1$ and $e_3=w_2w_3w_4=(88\pm4\sqrt2\,\I)/27$, while $e_1=2(1+\lambda)=4$ and $e_2=6\lambda=6$.
By complex conjugation, we may assume that $e_3=(88+4\sqrt2\,\I)/27$.
Then
\begin{equation*}
 F(w)=w\biggl(w-\frac{4+\sqrt2\,\I}{3}\biggr)^2
          \biggl(w-\frac{4-2\sqrt2\,\I}{3}\biggr).
\end{equation*}
Since $\abs{\zeta_1-a}=2/\sqrt3$, the identity $f(z)=(\zeta_1-a)^4F((z-a)/(\zeta_1-a))$ gives, for some $b\in\C$ and $\theta\in\R$, \begin{equation*}f(z)=\bigl(z-b-e^{\I\theta}\bigr)^2 \biggl((z-b)^2+\frac23e^{\I\theta}(z-b)+e^{2\I\theta}\biggr).\end{equation*}
Its zeros $b+e^{\I\theta}$ and $b+e^{\I\theta}\bigl(-\frac13\pm\frac{2\sqrt2}3\I\bigr)$ lie in $\D$, so direct expansion gives \begin{equation*}1+\abs b^2=\frac14\Absbig{b+e^{\I\theta}}^2 +\frac38\Absbigg{b+e^{\I\theta}\biggl(-\frac13+\frac{2\sqrt2}3\I\biggr)}^2 +\frac38\Absbigg{b+e^{\I\theta}\biggl(-\frac13-\frac{2\sqrt2}3\I\biggr)}^2\le1.\end{equation*}
Thus $b=0$, and $f(z)=\bigl(z-e^{\I\theta}\bigr)^2\bigl(z^2+\frac23e^{\I\theta}z+e^{2\I\theta}\bigr)$.
For a polynomial $f$ with leading coefficient $c\ne0$, the preceding argument applies to $f/c$ and shows that $f$ has the form~\eqref{eq:quartic-extremal}.
\end{proof}
\smallskip

\section{Proof of Theorem~\ref{thm:completeweighted}}\label{sec:generalproof}

Fix integers $s\ge2$ and $r\ge1$, and let $g$, $a$, and $G$ be as in Theorem~\ref{thm:completeweighted}. We may assume that $g$ is monic, $a\in[0,1]$, and $A(G)\subseteq\D$.
Indeed, the zeros of $g$ and $G$ are unchanged by normalization, and hence $g$ may be assumed to be monic. Write $a=\abs{a}e^{\I\theta}$, with $\theta=0$ if $a=0$, and set $\widetilde g(w)\=e^{-\I r\theta}g\bigl(e^{\I\theta}w\bigr)$. Then $\widetilde g$ is monic and $(w-\abs{a})\widetilde g'(w)+s\widetilde g(w)=e^{-\I r\theta}G\bigl(e^{\I\theta}w\bigr)$. Consequently, we may assume that $a\in[0,1]$, since this rotation preserves the distances from $a$ to the zeros of $G$. Finally, since $((z-a)^s g(z))'=(z-a)^{s-1}G(z)$, the Gauss--Lucas theorem implies $A(G)\subseteq\D$.

\smallskip
For ease of reference, we recall and fix the notation used throughout Section~\ref{sec:generalproof}, unless otherwise stated:
\begin{equation*}
 s,\,r\in\N,\quad s\ge2,\quad a\in[0,1],\quad
 N\=r+s,\quad R_s\=\frac{2s}{s+1}.
\end{equation*}
Write
\begin{equation*}
 g(z)\=\prod_{j=1}^r(z-\eta_j),\quad \text{ where }\eta_j\in\D,
\end{equation*}
and
\begin{equation*}
 G(z)\=(z-a)g'(z)+sg(z)=N\prod_{j=1}^r(z-\zeta_j),\quad\text{ where }\zeta_j\in\D.
\end{equation*}
\smallskip

Then, the conclusion of Theorem~\ref{thm:completeweighted} is as follows:
\begin{equation*}
 \min_{1\le j\le r}\abs{a-\zeta_j}\le R_s.
\end{equation*}

\subsection{Integral identities and estimates}\label{sec:identities}

We first establish a necessary condition for the existence of a counterexample.

\begin{lemma}\label{lem:reduction}
Suppose $\abs{a-\zeta_j}>R_s$ for $j=1,\,2,\,\ldots,\,r$.
Then $R_s-1<a<1$ and $g(a)\ne0$.
\end{lemma}
\begin{proof}
Since $\abs{a-\zeta_j}>R_s$ for $j=1,\,2,\,\ldots,\,r$, we have $G(a)=sg(a)\ne0$ and hence $g(a)\ne0$. Moreover, $R_s<\abs{a-\zeta_j}\le a+1$ for $j=1,\,2,\,\ldots,\,r$, so $a>R_s-1$.

If $a=1$ and $g(1)\ne0$,
\begin{equation*}
 \sum_{j=1}^r\frac1{1-\zeta_j}
 =\frac{G'(1)}{G(1)}
 =\frac{s+1}{s}\frac{g'(1)}{g(1)}
 =\frac{s+1}{s}\sum_{j=1}^r\frac1{1-\eta_j}.
\end{equation*}
Thus,
\begin{equation*}
 \sum_{j=1}^r\Rea\frac1{1-\zeta_j}
  =\frac{s+1}{s}\sum_{j=1}^r\Rea\frac1{1-\eta_j}
  \ge\frac{(s+1)r}{2s},
\end{equation*}
where we used a basic fact that $\Rea(1/(1-\eta_j))\ge1/2$ for $\abs{\eta_j}\le1$, $\eta_j\ne1$, and $j=1,\,2,\,\ldots,\,r$.
Then $\Rea(1/(1-\zeta_{j_0}))\ge 1/R_s$ for some $j_{0}\in\{1,\,2,\,\ldots,\,r\}$, and hence $1/\abs{1-\zeta_{j_0}}\ge1/R_s$, a contradiction.

This proves the lemma.
\end{proof}

\smallskip

Throughout the remainder of Section~\ref{sec:generalproof}, unless otherwise stated, we assume that $\abs{a-\zeta_j}>R_s$ for $j=1,\,2,\,\ldots,\,r$ and use the following notation.
\begin{equation}\label{eq:notation}
 k\=\frac{s+1}{2s},\quad q_j\=\frac1{a-\zeta_j},\quad
 x\=\frac1r\sum_{j=1}^r\Rea q_j,\quad
 \delta\=1-a^2,\quad \alpha\=\frac{r\delta}{2}.
\end{equation}
The hypothesis and Lemma~\ref{lem:reduction} imply
\begin{equation}\label{eq:counterexamplebounds}
 0<a<1,\quad \delta>0,\quad \alpha>0,\quad
 \abs{q_j}=\frac1{\abs{a-\zeta_j}}<k,
\end{equation}
where $j=1,\,2,\,\ldots,\,r$.
Since $G(a)\ne0$, for each $v\in\C$,
\begin{equation*}
\begin{aligned}
 \frac{G(a+v)}{G(a)}
 &=\prod_{j=1}^r\frac{a+v-\zeta_j}{a-\zeta_j}
 =\prod_{j=1}^r(1+vq_j).
\end{aligned}
\end{equation*}
For $t\in\C$ and $m\in \mathbb{N}_0$, define
\begin{equation*}
 F(t)\=\frac{G(a-at)}{G(a)}
       =\prod_{j=1}^r(1-atq_j)\quad\text{ and }\quad
 I_m\=\int_0^1\! t^m F(t)\,\mathrm{d}t.
\end{equation*}

We adopt Tao's product representation for $F$, while $I_m$ extends the integral in his first origin identity by introducing the weight $t^m$. We also obtain the following lemma, which gives weighted versions of his first origin and polar identities~\cite[Lemma~6]{Tao2026} with weight $t^{s-1}$ and provides the integral inequalities that will be used to exclude counterexamples.

\begin{lemma}\label{lem:identities}
Suppose $\abs{a-\zeta_j}>R_s$ for $j=1,\,2,\,\ldots,\,r$.
Then for each $u\in\C$,
\begin{equation}\label{eq:integral}
\begin{aligned}
 \frac{g(a-au)}{g(a)}
 &=s\int_0^1\!t^{s-1}F(ut)\,\mathrm{d}t.
\end{aligned}
\end{equation}
In particular,
\begin{equation}\label{eq:values}
 g(0)=\frac{N I_{s-1}}{\prod_{j=1}^r q_j}\quad\text{ and }\quad
 g(1/a)=\frac{s g(a)}{a^r}\int_0^1\!t^{s-1}\prod_{j=1}^r(a+\delta t q_j)\,\mathrm{d}t.
\end{equation}
Consequently,
\begin{equation}\label{eq:integralbounds}
 \abs{I_{s-1}}\le\frac{k^r}{N}\quad\text{ and }\quad
 1\le s\int_0^1\!t^{s-1}
       \bigl(a^2+2a\delta xt+\delta^2 k^2t^2\bigr)^{r/2}\,\mathrm{d}t.
\end{equation}
\end{lemma}

\begin{proof}
Since $(u^s g(a+u))^{\prime}=u^{s-1}G(a+u)$, for each nonzero $u\in \C$,
\begin{equation*}
 u^s g(a+u)=\int_0^u\! v^{s-1}G(a+v)\,\mathrm{d}v
           =u^s\int_0^1\!t^{s-1}G(a+ut)\,\mathrm{d}t.
\end{equation*}
Thus, for each $u\in\C$,
\begin{equation*}
 g(a+u)=\int_0^1\!t^{s-1}G(a+ut)\,\mathrm{d}t.
\end{equation*}
The equality at $u=0$ follows from the definition of $G(z)$. Replacing $u$ by $-au$ and recalling that $G(a-aut)=sg(a)F(ut)$, we obtain \eqref{eq:integral}.

Substituting $u\=1$ in \eqref{eq:integral} gives the first identity in \eqref{eq:values}.
Since $\abs{g(0)}=\prod_{j=1}^{r}\abs{\eta_j}\le1$ and $\abs{q_j}<k$ for $j=1,\,2,\,\ldots,\,r$, the first inequality in \eqref{eq:integralbounds} follows.

Since $a>0$ by~\eqref{eq:counterexamplebounds} and $a+\frac{\delta}{a}=\frac{1}{a}$, substituting $u=-\frac{\delta}{a^2}$ in \eqref{eq:integral} gives
\begin{equation*}
 g(1/a)=s g(a)\int_0^1\!t^{s-1}F\biggl(-\frac{\delta t }{a^2} \biggr)\,\mathrm{d}t
 =\frac{s g(a)}{a^r}\int_0^1\!t^{s-1}\prod_{j=1}^r(a+\delta t q_j)\,\mathrm{d}t,
\end{equation*}
which proves the second identity in \eqref{eq:values}. By the definition of $\eta_j$,
\begin{equation*}
 g(1/a)=\frac{g(a)}{a^r}\prod_{j=1}^r\frac{1-a\eta_j}{a-\eta_j}.
\end{equation*}
For $j=1,\,2,\,\ldots,\,r$, we have $\abs{1-a\eta_j}^2-\abs{a-\eta_j}^2 =\bigl(1-a^2\bigr)\bigl(1-\abs{\eta_j}^2\bigr)\ge0$.
Thus,
\begin{equation*}
 1\le\Absbigg{\frac{a^r g(1/a)}{g(a)}}
 \le s\int_0^1\!t^{s-1}\prod_{j=1}^r\abs{a+\delta t q_j}\,\mathrm{d}t.
\end{equation*}
To estimate the product in this integral, we use the arithmetic--geometric mean inequality (abbreviated as AM--GM).
In particular, for $t\in[0,1]$,
\begin{equation*}
 \prod_{j=1}^r\abs{a+\delta t q_j}
 \le\biggl(\frac1r\sum_{j=1}^r\abs{a+\delta t q_j}^2\biggr)^{r/2}.
\end{equation*}
The definition of $x$ and the bound $\abs{q_j}<k$ give
\begin{equation*}
 \frac1r\sum_{j=1}^r\abs{a+\delta t q_j}^2
 =a^2+2a\delta tx+\frac{\delta^2t^2}{r}\sum_{j=1}^r\abs{q_j}^2
 \le a^2+2a\delta tx+\delta^2k^2t^2,
\end{equation*}
which proves the second inequality in \eqref{eq:integralbounds}.
\end{proof}
\smallskip

To estimate $\abs{F(t)}$ on $[0,1]$, we again apply AM--GM to the product defining $F(t)$. For $0\le t\le1$, define
\begin{equation}\label{eq:betadefinition}
 \beta(t)\=1-2axt+a^2k^2t^2.
\end{equation}
The bounds $\abs{q_j}<k$ in~\eqref{eq:counterexamplebounds} imply
\begin{equation}\label{eq:betaestimation}
 \frac1r\sum_{j=1}^r\abs{1-atq_j}^2=1-2axt+\frac{a^2t^2}{r}\sum_{j=1}^r\abs{q_j}^2\le\beta(t).
\end{equation}
Thus, AM--GM gives $\abs{F(t)}\le\beta(t)^{r/2}$ and the following lemma shows further bounds on $x$ and $\beta$.

\begin{lemma}\label{lem:disk}
Suppose $\abs{a-\zeta_j}>R_s$ for $j=1,\,2,\,\ldots,\,r$. Then
\begin{equation}\label{eq:xbounds}
 \bigl( 1-\delta k^2 \bigr) \big/(2a)\le x\le k.
\end{equation}
For $0\le t\le1$,
\begin{equation}\label{eq:diskbounds}
 0\le\beta(t)\le1-\bigl(1-k^2\bigr)t-a^2 k^2t(1-t)\le1.
\end{equation}
If $s\ge3$, one also has
\begin{equation}\label{eq:betau}
 \beta(t)\le1-at+a^2k^2t^2.
\end{equation}
\end{lemma}
\begin{proof}
The upper bound in \eqref{eq:xbounds} follows by the definition of $x$ and $\Rea q_j\le\abs{q_j}<k$. Since $aq_j-1=\zeta_j q_j$ and $\abs{\zeta_j}\le1$, we have $a^2\abs{q_j}^2-2a\Rea q_j+1 =\abs{aq_j-1}^2\le\abs{q_j}^2$, and thus $ 2a\Rea q_j\ge1-\delta\abs{q_j}^2$.
Hence $\abs{q_j}<k$ yields
\begin{equation*}
 2ax=\frac{2a}{r}\sum_{j=1}^r\Rea q_j
 \ge1-\frac{\delta}{r}\sum_{j=1}^r\abs{q_j}^2
 \ge1-\delta k^2.
\end{equation*}
Then the lower bound in \eqref{eq:xbounds} follows by noticing that $a>0$.

For $0\le t\le1$, the upper bound for $x$ gives $\beta(t)\ge(1-akt)^2\ge0$, while the lower bound for $2ax$ gives
$ \beta(t)\le1- \bigl(1-\delta k^2\bigr)t+a^2k^2t^2
 =1-\bigl(1-k^2\bigr)t-a^2k^2t(1-t)\le1$.
This proves \eqref{eq:diskbounds}.

Finally, suppose $s\ge3$. Then $k^2=\bigl(\frac{s+1}{2s}\bigr)^2\le\frac12$.
Since $0<a<1$ from Lemma \ref{lem:reduction},
\begin{equation*}
 \delta k^2
 =(1-a)(1+a)k^2
 \le (1-a)(1+a)/2
 \le1-a.
\end{equation*}
Together with \eqref{eq:xbounds}, this gives $2ax\ge1-\delta k^2\ge a$.
For $0\le t\le1$, the definition of $\beta$ therefore implies
$ \beta(t)=1-2axt+a^2k^2t^2
 \le1-at+a^2k^2t^2$,
which proves \eqref{eq:betau}.
\end{proof}
\smallskip

\subsection{Four cases}\label{sec:cases}
We now prove Theorem~\ref{thm:completeweighted} by partitioning the integer pairs $(s,r)$ with $s\ge2$ and $r\ge1$ into four disjoint cases. In each case, we will assume that there exists a counterexample and derive a contradiction. The parameter ranges and the corresponding propositions are listed in Table~\ref{tab:cases}.
\begin{table}[ht]
\centering
\begin{tabular}{ccc}
\toprule
case&parameter range&result\\
\midrule
I&$r\le(s+1)^2$&Proposition~\ref{prop:small}\\
II&$r>(s+1)^2,\ s\ge9$&Proposition~\ref{prop:large}\\
III&$2\le s\le8,\ r\ge r_0(s)$&Proposition~\ref{prop:tail}\\
IV&$2\le s\le8,\ (s+1)^2<r<r_0(s)$&Proposition~\ref{prop:finiteanalytic}\\
\bottomrule
\end{tabular}
\par\vspace{8pt}
\small
\caption[Parameter ranges]{$r_0(s)$ is the degree threshold
separating Cases~III and~IV, specified in Table~\ref{tab:tail}.}
\label{tab:cases}
\end{table}
\FloatBarrier

Under the counterexample hypothesis, define
\begin{equation*}
 \overline A\=ark\quad \text{ and }\quad \Sfun(u)\=\sum_{j=0}^{s-1}\frac{u^j}{j!},
\end{equation*}
for $u\in\C$. Set
\begin{equation}\label{eq:U}
\begin{aligned}
 M\=\Sfun(\overline A)\beta(1)^{r/2}
   +\frac{\overline A^s k^r}{(s-1)!(r+s)},\quad
 E\=a^2rk^2\int_0^1\!t\Sfun(\overline A t)
          \beta(t)^{(r-1)/2}\,\mathrm{d}t,\quad
 U\=M+E.
\end{aligned}
\end{equation}

We first prove the following lower bound for $U$, which is a necessary condition for the existence of a counterexample. This bound will be used to derive contradictions in the subsequent arguments.

\begin{lemma}\label{lem:originineq}
Suppose $\abs{a-\zeta_j}>R_s$ for $j=1,\,2,\,\ldots,\,r$.
Then $U\ge 1$.
\end{lemma}
\begin{proof}
Write $A_0\=a\sum_{j=1}^r q_j$ and $F_j(t)\=\prod_{\substack{1\le\ell\le r\\ \ell\ne j}}(1-atq_\ell)$ for $j=1,\,2,\,\ldots,\,r$ and $t\in[0, 1]$.
Since $F(t)=(1-atq_j)F_j(t)$, we have $F_j(t)=F(t)+atq_jF_j(t)$. Then
\begin{equation}\label{eq:Fderivative}
 F'(t)=-a\sum_{j=1}^r q_jF_j(t)
      =-a\sum_{j=1}^r q_j(F(t)+atq_jF_j(t))
      =-A_0F(t)-a^2t\sum_{j=1}^r q_j^2F_j(t).
\end{equation}
For each $m\in\N_0$, let $E_m\=a^2\int_0^1\!t^{m+1}\sum_{j=1}^r q_j^2F_j(t)\,\mathrm{d}t$.
Then \eqref{eq:Fderivative} and the definition of $I_m$ give that
\begin{equation*}
 \int_0^1\!t^mF'(t)\,\mathrm{d}t=-A_0I_m-E_m.
\end{equation*}
When $m=0$, the left-hand side is $F(1)-F(0)=F(1)-1$. For $m\ge1$,
\begin{equation*}
 \int_0^1\!t^mF'(t)\,\mathrm{d}t
 =\left.t^mF(t)\right|_0^1
   -m\int_0^1\!t^{m-1}F(t)\,\mathrm{d}t
 =F(1)-mI_{m-1}.
\end{equation*}
Thus, we have
\begin{equation}\label{eq:Imrecurrence}
 1=F(1)+A_0 I_0+E_0\quad\text{ and }\quad
 m I_{m-1}=F(1)+A_0 I_m+E_m,
\end{equation}
for $m=0$ and $m\ge1$, respectively.
Multiplying the second equality in \eqref{eq:Imrecurrence} by $A_0^m/m!$, summing over $1\le m\le s-1$, and adding the first equality, we obtain
\begin{equation*}
 1+\sum_{m=1}^{s-1}\frac{A_0^m}{(m-1)!}I_{m-1}
 =\Sfun(A_0)F(1)
  +\sum_{m=0}^{s-1}\frac{A_0^{m+1}}{m!}I_m
  +\sum_{m=0}^{s-1}\frac{A_0^m}{m!}E_m.
\end{equation*}
The sum on the left is equal to $\sum_{m=0}^{s-2}A_0^{m+1}I_m/m!$, so it cancels all but the last term of the sum involving $I_m$ on the right. Hence,
\begin{equation}\label{eq:telescoping}
 1=\Sfun(A_0)F(1)+\frac{A_0^s}{(s-1)!}I_{s-1}
      +\sum_{m=0}^{s-1}\frac{A_0^m}{m!}E_m.
\end{equation}

To estimate $\sum_{j=1}^r\abs{F_j(t)}$, we first prove the following auxiliary inequalities. For nonnegative numbers $b_1,\,b_2,\,\ldots,\,b_r$, we have
\begin{equation*}
 \sum_{j=1}^r\prod_{\substack{1\le\ell\le r\\ \ell\ne j}}b_\ell
 \le r\biggl(\frac1r\sum_{j=1}^r b_j\biggr)^{r-1}
 \le r\biggl(\frac1r\sum_{j=1}^r b_j^2\biggr)^{(r-1)/2}.
\end{equation*}
Indeed, for $r\ge2$, the two inequalities follow from Maclaurin's inequality and the Cauchy--Schwarz inequality, respectively. For $r=1$, the left-hand side equals $1$ by the empty product convention.

For $0\le t\le1$, applying both inequalities above with $b_j\=\abs{1-atq_j}$ and using $\frac{1}{r}\sum_{j=1}^r b_j^2\le\beta(t)$ from \eqref{eq:betaestimation}, we get
\begin{equation}\label{eq:omittedproductbounds}
 \sum_{j=1}^r\abs{F_j(t)}\le r\beta(t)^{(r-1)/2}.
\end{equation}
Substituting $t=1$ in the AM--GM bound $\abs{F(t)}\le\beta(t)^{r/2}$ gives
\begin{equation}\label{eq:Fendpointbound}
 \abs{F(1)}\le\beta(1)^{r/2}.
\end{equation}
Since $\abs{A_0}\le a\sum_{j=1}^r\abs{q_j}\le\overline A$, we have $\abs{\Sfun(A_0)}\le\Sfun(\overline A)$. It follows from $\abs{q_j}<k$ and \eqref{eq:omittedproductbounds} that
\begin{equation}\label{eq:Embound}
 \Absbigg{\sum_{m=0}^{s-1}\frac{A_0^m}{m!}E_m}
 \le a^2\int_0^1\!t\Sfun(\abs{A_0}t)
          \sum_{j=1}^r\abs{q_j}^2\abs{F_j(t)}\,\mathrm{d}t
 \le a^2rk^2\int_0^1\!t\Sfun(\overline A t)
          \beta(t)^{(r-1)/2}\,\mathrm{d}t=E.
\end{equation}
By combining the first inequality in \eqref{eq:integralbounds}, \eqref{eq:telescoping}, \eqref{eq:Fendpointbound}, and \eqref{eq:Embound}, we conclude
\begin{equation*}
\begin{aligned}
 1\le\abs{\Sfun(A_0)}\abs{F(1)}
      +\frac{\abs{A_0}^s}{(s-1)!}\abs{I_{s-1}}+E
  \le\Sfun(\overline A)\beta(1)^{r/2}
      +\frac{\overline A^s k^r}{(s-1)!(r+s)}+E
   =M+E=U,
\end{aligned}
\end{equation*}
which proves the lemma.
\end{proof}

\smallskip

We next derive an upper bound for $M$ that will be used in Cases~II--IV. For $r>(s+1)^2$, we have $rk\ge1$.
The estimate $\beta(1)\le k^2$ from \eqref{eq:diskbounds} and $a\le1$ give
\begin{equation}\label{eq:Mbound}
 M\le k^r\biggl(\Sfun(rk)+\frac{(rk)^s}{(s-1)!(r+s)}\biggr)
 \le (s+k)(rk)^{s-1}k^r,
\end{equation}
where we used $\Sfun(rk)\le s(rk)^{s-1}$ and $rk/((s-1)!(r+s))\le k$ in the last inequality.

For each fixed integer $s\ge2$ and each $j=0,\,1,\,\ldots,\,s$, the sequence $r^jk^r$ is strictly decreasing for integers $r\ge5s$. Indeed, $k=(s+1)/(2s)\le3/4$ and
\begin{equation}\label{eq:degreeMonotonicity}
 \frac{(r+1)^jk^{r+1}}{r^jk^r}
 =k\biggl(1+\frac1r\biggr)^j
 =k\sum_{\ell=0}^{j}\binom{j}{\ell}\frac1{r^\ell}
 \le k\sum_{\ell=0}^{j}\Bigl(\frac sr\Bigr)^\ell
 \le k\sum_{\ell=0}^{+\infty}\Bigl(\frac sr\Bigr)^\ell
 =\frac{k}{1-\frac{s}{r}}\le\frac{15}{16}<1.
\end{equation}
In particular, the upper bound $(s+k)(rk)^{s-1}k^r$ in~\eqref{eq:Mbound} is strictly decreasing for integers $r\ge5s$.

The estimate for $E$ depends on the decay of $\beta(t)$.
The next lemma improves the lower bound for $x$, hence strengthens the linear term in $\beta(t)$ without changing its quadratic coefficient.
This improvement will be used in Case~IV.

In the sequel, $\E$ denotes expectation with respect to the distribution specified in each statement.

\begin{lemma}\label{lem:mean}
Suppose $\abs{a-\zeta_j}>R_s$ for $j=1,\,2,\,\ldots,\,r$. Then
\begin{equation}\label{eq:mean}
 2ax>1+\frac1{s+\alpha}-\delta k^2.
\end{equation}
\end{lemma}
\begin{proof}
Let $C\=2ax+\delta k^2$.
To prove \eqref{eq:mean}, it suffices to show $C>1+1/(s+\alpha)$.
Let $T$ be a random variable taking values in $[0,1]$ with probability density $s t^{s-1}$.
The second inequality in \eqref{eq:integralbounds} then reads
\begin{equation*}
 1\le\E\bigl[\bigl(a^2+2a\delta xT+\delta^2k^2T^2\bigr)^{r/2}\bigr].
\end{equation*}
For $0\le t\le1$, the definitions of $\delta$ and $\alpha$, together with $t^2\le t$, give
\begin{equation*}
 (r/2)\bigl(a^2+2a\delta xt+\delta^2k^2t^2-1\bigr)
 =\alpha\bigl(-1+2axt+\delta k^2t^2\bigr)
 \le\alpha(-1+Ct).
\end{equation*}
Applying $w^{r/2}\le e^{(r/2)(w-1)}$ for $w\ge0$ to the integrand therefore yields
\begin{equation}\label{eq:meanexp}
 1\le\E \bigl[e^{\alpha(-1+CT)}\bigr].
\end{equation}
To bound this expectation, set $U_0\=1-T$, so $-1+CT=(C-1)-CU_0$. Then
\begin{equation}\label{eq:U0mean}
 \E [U_0]=s\int_0^1\!(1-t)t^{s-1}\,\mathrm{d}t
 =\frac1{s+1}.
\end{equation}
Write $y_+\=\max\{y,\,0\}$ for real $y$. For $0\le u\le1$,
\begin{equation}\label{eq:U0positivepart}
 \E[(U_0-u)_+]
 =s\int_0^{1-u}\!(1-u-t)t^{s-1}\,\mathrm{d}t
 =\frac{(1-u)^{s+1}}{s+1}
 \le\frac{e^{-(s+1)u}}{s+1},
\end{equation}
where we used a basic inequality $1-u\le e^{-u}$ in the last inequality. Since $U_0\le1$, we have $\E[(U_0-u)_+]=0$ for $u\ge1$. Thus the upper bound in \eqref{eq:U0positivepart} holds for every $u\ge0$. Besides, for $b\ge0$ and $y\ge0$, we have
$ e^{-by}=1-by+\int_0^{+\infty}\! b^2 e^{-bu}(y-u)_+\,\mathrm{d}u$.
Taking expectations in this identity, and applying \eqref{eq:U0mean} and \eqref{eq:U0positivepart}, we obtain
\begin{equation}\label{eq:laplace}
\begin{aligned}
 \E \bigl[e^{-bU_0}\bigr]
 &=1-b\E [U_0]+\int_0^{+\infty}\!b^2e^{-bu}\E[(U_0-u)_+]\,\mathrm{d}u\\
 &\le1-\frac{b}{s+1}+\frac{b^2}{(s+1)(s+1+b)}
 =\frac{s+1}{s+1+b}.
\end{aligned}
\end{equation}
Suppose, for a contradiction, that $C\le1+1/(s+\alpha)$. Let $v\=\alpha/(s+\alpha)$.
Since $\alpha>0$ by~\eqref{eq:counterexamplebounds}, we have $0<v<1$ and $\alpha C\le\alpha+v$. It follows from $T\ge0$, $T=1-U_0$, and \eqref{eq:laplace} with $b=\alpha+v>0$ that
\begin{equation}\label{eq:Cthresholdbound}
 \E \bigl[e^{\alpha(-1+CT)}\bigr]
 \le\E \bigl[e^{-\alpha+(\alpha+v)T}\bigr]
 =e^v\E \bigl[e^{-(\alpha+v)U_0}\big]
 \le\frac{e^v(s+1)}{s+1+\alpha+v}
 =\frac{e^v(1-v)(s+1)}{s+1-v^2}.
\end{equation}
We claim that $\frac{e^v(1-v)(s+1)}{s+1-v^2}<1$, which, together with~\eqref{eq:Cthresholdbound}, contradicts~\eqref{eq:meanexp}. Hence the assertion follows.

Indeed, $1-\frac{e^v(1-v)(s+1)}{s+1-v^2} =\frac{s+1}{s+1-v^2}\bigl(1-\frac{v^2}{s+1}-e^v(1-v)\bigr)$. The first factor is positive. The second vanishes at $v=0$ and has derivative $\frac{\mathrm{d}}{\mathrm{d}v}\bigl(1-\frac{v^2}{s+1}-e^v(1-v)\bigr) =v\bigl(e^v-\frac{2}{s+1}\bigr)>0$ for $s\ge2$ and $0<v<1$. This proves the claim.
\end{proof}

\subsubsection{Case I}\label{sec:small}
Assume $r\le(s+1)^2$. Under the counterexample hypothesis, by applying $x\le k$ in the second inequality of \eqref{eq:integralbounds}, we have $1\le s\int_0^1\! t^{s-1}\bigl(a+k\bigl(1-a^2\bigr)t\bigr)^r\,\mathrm{d}t$. The following lemma will show that this integral is strictly less than $1$, which then excludes Case~I.

\begin{lemma}\label{lem:moment}
Let $s\in\N$, $p,\,A,\,B\in\R$ satisfy $p\ge3$ and $A\ge B\ge0$. Let $T$ be a random variable taking values in $[0,1]$ with probability density $s t^{s-1}$ and $Y\=(s+1)\frac{T}{s}-1$. Then
\begin{equation}\label{eq:moment}
 \E[(A+BY)^p]^{2/p}\le A^2+\frac{p-1}{s(s+2)}B^2.
\end{equation}
\end{lemma}
\begin{proof}
If $A=0$, then $B=0$ and the result is immediate. Now we assume $A>0$. Since the density of $Y$ on $[-1,1/s]$ is $p_Y(y)=\frac{s^{s+1}}{(s+1)^s}(1+y)^{s-1}$, a direct calculation gives $\E [Y]=0$ and $\E \bigl[Y^2\bigr]=1/(s(s+2))$.
Define $\tau(y)\=(1+y)(1-sy)/(s(s+1))$ for $y\in[-1,1/s]$. Then
\begin{equation}\label{eq:tauendpoints}
 \tau(y)\ge0\quad\text{ and }\quad (\tau p_Y)(-1)=(\tau p_Y)(1/s)=0.
\end{equation}
Moreover, $\E[\tau(Y)]=1/(s(s+2))$, $(\tau p_Y)'(y)=-y p_Y(y)$, and
\begin{equation}\label{eq:comparisonderivative}
 ((sy+s-1)\tau(y)p_Y(y))'
 =\left(\frac{\tau(y)}{\E[\tau(Y)]}-1\right)p_Y(y).
\end{equation}

For $0\le b\le B$, let $h(y)\=(A+by)^{p-2}$.
We claim that
\begin{equation}\label{eq:weightedcomparison}
 \E[\tau(Y)h(Y)]\le\E[\tau(Y)]\,\E [h(Y)].
\end{equation}
Since $p\ge3$ and $A+by\ge A-b\ge A-B\ge0$ on $[-1,1/s]$, $h$ is continuously differentiable and $h'$ is nondecreasing.
Hence $(sy+s-1)(h'(y)-h'(1/s-1))\ge0$ for $y\in[-1,1/s]$.
Multiplying by $-\tau(y)\le0$ and taking expectations gives
\begin{equation}\label{eq:expectationcomparison}
 -\E[\tau(Y)(sY+s-1)h'(Y)]
 \le-h'(1/s-1)\E[\tau(Y)(sY+s-1)].
\end{equation}
Since $h'\ge0$ and $s\ge1$, it follows from \eqref{eq:tauendpoints}, \eqref{eq:comparisonderivative}, and \eqref{eq:expectationcomparison} that
\begin{align*}
 \frac{\E[\tau(Y)h(Y)]}{\E[\tau(Y)]}-\E [h(Y)]
 &=\int_{-1}^{1/s}\! h(y)\biggl(\frac{\tau(y)}{\E[\tau(Y)]}-1\biggr)p_Y(y)\,\mathrm{d}y\notag\\
 &=\int_{-1}^{1/s}\! h(y)((sy+s-1)\tau(y)p_Y(y))'\,\mathrm{d}y\notag\\
 &=-\E[\tau(Y)(sY+s-1)h'(Y)]\notag\displaybreak[1]\\
 &\le-h'(1/s-1)\E[\tau(Y)(sY+s-1)]\\
 &=-h'(1/s-1)\frac{s^s}{(s+1)^{s+1}}
 \int_{-1}^{1/s}\! (sy+s-1)(1+y)^s(1-sy)\,\mathrm{d}y\notag\\
 &=-h'(1/s-1)\frac{s-1}{s(s+3)}\notag\\
 &\le0.\notag
\end{align*}
Thus \eqref{eq:weightedcomparison} holds.

Let $\Psi(b)\=\E[(A+bY)^p]$ for $0\le b\le B$.
By combining $(\tau p_Y)'(y)=-y p_Y(y)$, \eqref{eq:tauendpoints}, \eqref{eq:weightedcomparison}, and H\"older's inequality with conjugate exponents $p/(p-2)$ and $p/2$, we obtain

\begin{align*}
 \Psi'(b)
 &=p\E\bigl[Y(A+bY)^{p-1}\bigr]\\
 &=-p\int_{-1}^{1/s}(A+by)^{p-1}\! (\tau p_Y)'(y)\,\mathrm{d}y\\
 &=p(p-1)b\E\bigl[\tau(Y)(A+bY)^{p-2}\bigr]\\
 &\le p(p-1)b\,\E[\tau(Y)]\,\E\bigl[(A+bY)^{p-2}\bigr]\\
 &\le p(p-1)b\,\E[\tau(Y)]\,
 \E[(A+bY)^p]^{(p-2)/p}\E[1]^{2/p}\\
 &=p(p-1)b\,\E[\tau(Y)]\,\Psi(b)^{(p-2)/p}.
\end{align*}
Since $A+by\ge A>0$ on $[0,1/s]$, $\Psi(b)>0$. Thus $\frac{\mathrm{d}}{\mathrm{d}b}\Psi(b)^{2/p} =\frac{2}{p}\Psi(b)^{(2-p)/p}\Psi'(b) \le2(p-1)b\,\E[\tau(Y)]$. Since $\Psi(0)=A^p$ and $\E[\tau(Y)]=1/(s(s+2))$,
\begin{equation*}
 \int_0^B\!\frac{\mathrm{d}}{\mathrm{d}b}\Psi(b)^{2/p}\,\mathrm{d}b
 =\Psi(B)^{2/p}-A^2
 \le2(p-1)\E[\tau(Y)]\int_0^B\! b\,\mathrm{d}b
 =\frac{p-1}{s(s+2)}B^2,
\end{equation*}
which proves \eqref{eq:moment}.
\end{proof}
\smallskip

\begin{proposition}\label{prop:small}
Under the hypotheses of Theorem~\ref{thm:completeweighted}, suppose $1\le r\le(s+1)^2$. Then $G$ has a zero $\zeta$ satisfying $\abs{a-\zeta}\le R_s$.
\end{proposition}
\begin{proof}
Suppose, for a contradiction, that $\abs{a-\zeta_j}>R_s$ for $j=1,\,2,\,\ldots,\,r$. By combining the second inequalities in \eqref{eq:integralbounds} and \eqref{eq:xbounds}, we obtain
\begin{equation*}
 1\le\E\bigl[\bigl(a+k\bigl(1-a^2\bigr)T\bigr)^r\bigr].
\end{equation*}
Let $B\=\bigl(1-a^2\bigr)/2$ and $A\=a+B$. Then $a+k\bigl(1-a^2\bigr)T=A+BY$, $1-A^2-B^2=(1-a)^3(1+a)/2>0$, and $A\ge B\ge0$. H\"older's inequality and \eqref{eq:moment}, with $p=(s+1)^2$, give
\begin{equation*}
 \E\bigl[\bigl(a+k\bigl(1-a^2\bigr)T\bigr)^r\bigr]^{1/r}
 \le(\E[(A+BY)^p])^{1/p}
 \le\sqrt{A^2+B^2}<1,
\end{equation*}
which is a contradiction. Thus, the assertion follows.
\end{proof}

\subsubsection{Case II}\label{sec:large}
Assume $r>(s+1)^2$ and $s\ge9$. To estimate $E$, we first give the following bound for $\beta(t)$.

\begin{lemma}\label{lem:exponentiallarge}
Suppose $s\ge9$ and $\abs{a-\zeta_j}>R_s$ for $j=1,\,2,\,\ldots,\,r$. Then $\beta(t)\le e^{-at-(at)^2/10}$ for $0\le t\le1$.
\end{lemma}
\begin{proof}
Since $s\ge9$, we have $k=(s+1)/(2s)\le5/9$. Lemma~\ref{lem:disk} therefore implies $\beta(t)\le1-at+25(at)^2/81$. For $0\le t\le1$, let $f(t)\=\log\bigl(1-at+\frac{25}{81}(at)^2\bigr)+at+\frac{1}{10}(at)^2$. Then $f(0)=0$ and $f^{\prime}(t)=\frac{a^2t(-74+44at+25(at)^2)}{5(81-81at+25(at)^2)}\le0$, which proves the assertion.
\end{proof}

\smallskip

\begin{lemma}\label{lem:Elarge}
Suppose $s\ge9$, $r>(s+1)^2$, and $\abs{a-\zeta_j}>R_s$ for $j=1,\,2,\,\ldots,\,r$. Then $E<199/200$.
\end{lemma}
\begin{proof}
Let $b\=(r-1)/2$ and substitute $u\=at\ge 0$ in \eqref{eq:U}.
By Lemma~\ref{lem:exponentiallarge},
\begin{equation*}
\begin{aligned}
 E=rk^2\int_0^a\!u\Sfun(rku)\beta\Bigl(\frac{u}{a}\Bigr)^b\,\mathrm{d}u
  &\le rk^2\int_0^a\!u\Sfun(rku)e^{-bu-bu^2/10}\,\mathrm{d}u\\
  &\le rk^2\int_0^{+\infty}\!u\Sfun(rku)e^{-bu}
       \biggl(1-\frac{bu^2}{10}+\frac{b^2u^4}{200}\biggr)\,\mathrm{d}u.
\end{aligned}
\end{equation*}
The last inequality follows from $e^{-v}\le1-v+\frac{1}{2}v^2$ for $v\ge0$ and the positivity of the polynomial $1-\frac{1}{10}bu^2+ \frac{1}{200}b^2u^4$ for $u\ge0$.
By applying $\int_0^{+\infty}\!u^m e^{-bu}\,\mathrm{d}u=m!/b^{m+1}$ for any $m \in \mathbb{N}$, we obtain
\begin{equation}\label{eq:gamma}
\begin{aligned}
 E\le{}&rk^2\sum_{j=0}^{s-1}\frac{(rk)^j}{j!}
 \int_0^{+\infty}\!u^{j+1}e^{-bu}
 \biggl(1-\frac{bu^2}{10}+\frac{b^2u^4}{200}\biggr)\,\mathrm{d}u\\
 ={}&rk^2\sum_{j=0}^{s-1}\frac{(rk)^j}{j!}
 \biggl(\frac{(j+1)!}{b^{j+2}}-\frac{(j+3)!}{10b^{j+3}}
       +\frac{(j+5)!}{200b^{j+4}}\biggr)\\
 ={}&\frac{rk^2}{b^2}\sum_{j=0}^{s-1}(j+1)
 \biggl(\frac{rk}{b}\biggr)^j
 \biggl(1-\frac{(j+2)(j+3)}{10b}
  +\frac{(j+2)(j+3)(j+4)(j+5)}{200b^2}\biggr).
\end{aligned}
\end{equation}

Let $z\=rk/b>1$. To estimate the negative term in \eqref{eq:gamma}, we use Chebyshev's sum inequality: $\bigl(\sum_{j=0}^{s-1} w_jx_j\bigr)\bigl(\sum_{j=0}^{s-1} w_jy_j\bigr) \le\bigl(\sum_{j=0}^{s-1} w_j\bigr)\bigl(\sum_{j=0}^{s-1} w_jx_jy_j\bigr)$, for any nondecreasing real sequences $\{x_j\}_{j=0}^{s-1}$ and $\{y_j\}_{j=0}^{s-1}$ and positive sequence $\{w_j\}_{j=0}^{s-1}$~\cite{HardyLittlewoodPolya1952}.
Taking $w_j\=j+1$, $x_j\=z^j$, and $y_j\=(j+2)(j+3)$ gives
\begin{equation}\label{eq:negativecorrection}
\begin{aligned}
 -\frac1{10b}\sum_{j=0}^{s-1}(j+1)z^j(j+2)(j+3)
 &\le-\frac1{10b}
 \frac{\sum_{j=0}^{s-1}(j+1)(j+2)(j+3)}{\sum_{j=0}^{s-1}(j+1)}
 \sum_{j=0}^{s-1}(j+1)z^j\\
 &=-\frac{(s+2)(s+3)}{20b}\sum_{j=0}^{s-1}(j+1)z^j.
\end{aligned}
\end{equation}

Since $b=(r-1)/2$, $k=(s+1)/(2s)$, and $(s+1)/(r-1)<1$, we have
\begin{equation}\label{eq:geometricidentity}
\begin{aligned}
 \frac{rk^2}{b^2}\sum_{j=0}^{s-1}(j+1)\biggl(\frac{rk}{b}\biggr)^j
 &\le\frac{rk^2}{b^2}\Bigl(\frac{r}{r-1}\Bigr)^{s-1}
 \sum_{j=0}^{s-1}(j+1)(2k)^j\\
 &=\frac{(s+1)^2}{r}\Bigl(\frac{r}{r-1}\Bigr)^{s+1}
 =\frac{(s+1)^2}{r}\sum_{m=0}^{s+1}\binom{s+1}{m}\biggl(\frac{1}{r-1}\biggr)^{m}\\
 &\le\frac{(s+1)^2}{r}\sum_{m=0}^{+\infty}
 \biggl(\frac{s+1}{r-1}\biggr)^m
 =\frac{(s+1)^2(r-1)}{r(r-s-2)}.
\end{aligned}
\end{equation}

Combining \eqref{eq:gamma}, \eqref{eq:negativecorrection}, and \eqref{eq:geometricidentity} yields
\begin{equation*}
\begin{aligned}
 E&\le\frac{rk^2}{b^2}
 \biggl(1-\frac{(s+2)(s+3)}{10(r-1)}
 +\frac{(s+1)(s+2)(s+3)(s+4)}{50(r-1)^2}\biggr)
 \sum_{j=0}^{s-1}(j+1)\biggl(\frac{rk}{b}\biggr)^j\\
 &\le\frac{(s+1)^2(r-1)}{r(r-s-2)}
 \biggl(1-\frac{(s+2)(s+3)}{10(r-1)}
 +\frac{(s+1)(s+2)(s+3)(s+4)}{50(r-1)^2}\biggr).
\end{aligned}
\end{equation*}

To prove $E<199/200$, it suffices to show
\begin{equation*}
\begin{aligned}
 D(r)\={}&\frac{199}{200}(r-1)r(r-s-2)\\
       &-(s+1)^2\biggl((r-1)^2-\frac{(s+2)(s+3)(r-1)}{10}
       +\frac{(s+1)(s+2)(s+3)(s+4)}{50}\biggr)>0.
\end{aligned}
\end{equation*}
On the one hand,
\begin{equation*}
 D\bigl( (s+1)^2+1 \bigr)
 =\frac{(s+1)^3}{200} \bigl(15(s-9)^3+287(s-9)^2+1435(s-9)+427\bigr)>0.
\end{equation*}
On the other hand, since $(s+1)^2\ge4s$, we have $D''(r)\ge\frac{397(s+1)^2-199s}{100}>0$ and $D'\bigl( (s+1)^2+1 \bigr) \ge\frac{39s(s+1)^2}{20}-\frac{199(s+1)}{200}>0$. Thus $D$ is increasing when $r\ge (s+1)^2+1$, so the required bound follows.
\end{proof}

\smallskip

\begin{proposition}\label{prop:large}
Under the hypotheses of Theorem~\ref{thm:completeweighted}, suppose $s\ge9$ and $r>(s+1)^2$.
Then $G$ has a zero $\zeta$ satisfying $\abs{a-\zeta}\le R_s$.
\end{proposition}
\begin{proof}
Suppose, for a contradiction, that $\abs{a-\zeta_j}>R_s$ for $j=1,\,2,\,\ldots,\,r$.
It follows from \eqref{eq:Mbound}, \eqref{eq:degreeMonotonicity} with $j\=s-1$, and $k\=(s+1)/(2s)\le5/9$ that
\begin{align*}
 M\le(s+k)(rk)^{s-1}k^r
 &\le(s+k)\bigl(\bigl( (s+1)^2+1 \bigr)k\bigr)^{s-1}k^{(s+1)^2+1}\\
  &\le(s+1)(s+2)^{2s-2}(5/9)^{(s+1)^2+s}=:d_s.
\end{align*}
For $s\ge9$, we have $\bigl(1+\frac1{s+2}\bigr)^{2s-2}<e^2<8$ and
\begin{equation*}
 \frac{(s+4)^2(5/9)^{2s+6}}{(s+3)^2(5/9)^{2s+4}}
 =\biggl(1+\frac1{s+3}\biggr)^2\biggl(\frac59\biggr)^2
 \le\biggl(\frac{13}{12}\biggr)^2\biggl(\frac59\biggr)^2<1.
\end{equation*}
Thus, $(s+3)^2(5/9)^{2s+4}$ is strictly decreasing for integers $s\ge9$.
Since $(5/9)^{10}<2^{-8}$, we have
\begin{equation*}
\begin{aligned}
 \frac{d_{s+1}}{d_s}
 &=\frac{s+2}{s+1}(s+3)^2
   \biggl(1+\frac1{s+2}\biggr)^{2s-2}\biggl(\frac59\biggr)^{2s+4}\\
 &<\frac{11}{10}\,8(s+3)^2\biggl(\frac59\biggr)^{2s+4}
 \le\frac{11}{10}\,8\cdot12^2\biggl(\frac59\biggr)^{22}
 <\frac{6336}{5\cdot2^{16}}<1.
\end{aligned}
\end{equation*}
Thus,
$M\le d_s\le d_9=10\cdot11^{16}\bigl(\frac59\bigr)^{109}
 <\frac{10\cdot16^{16}}{2^{80}}=\frac{10}{2^{16}}<\frac1{200}$.
Then by applying Lemma~\ref{lem:Elarge}, we have $U=M+E<1$, which contradicts Lemma~\ref{lem:originineq}. Hence, the assertion follows.
\end{proof}

\subsubsection{Case III}\label{sec:tails}

For $2\le s\le8$, fix $r_0(s)$ and $c_s$ as in Table~\ref{tab:tail}. For fixed $s$, we also write $r_0\=r_0(s)$ for brevity. We also give the decay estimate that will be used to get upper bounds for $E$.

\begin{table}[ht]
\centering
\begin{tabular}{c*{7}{c}}
\toprule
$s$&2&3&4&5&6&7&8\\
\midrule
$r_0(s)$&150&120&110&72&98&128&162\\
$c_s$&$1/2$&$3/4$&$7/8$&1&1&1&1\\
\bottomrule
\end{tabular}
\par\vspace{8pt}
\small
\caption[Large-degree parameters]{$r_0(s)$ is the lower degree threshold for Case~III, and $c_s$ is the decay rate for the bound $\beta(t)\le e^{-c_s at}$ that will be proved in Lemma~\ref{lem:taildecay}.}\label{tab:tail}
\end{table}

\begin{lemma}\label{lem:taildecay}
Suppose $2\le s\le8$ and $\abs{a-\zeta_j}>R_s$ for $j=1,\,2,\,\ldots,\,r$.
Then $\beta(t)\le e^{-c_s at}$ for $0\le t\le1$, where $c_s$ is given in Table~\ref{tab:tail}.
\end{lemma}
\begin{proof}
Let $u\=at\in[0,1]$. For $s=2$, $c_2=1/2$, and $k^2=9/16$. Since $\frac7{16}+\frac9{16}a^2-\frac{31}{32}a =\frac9{16}\bigl(a-\frac{31}{36}\bigr)^2+\frac{47}{2304}>0$, \eqref{eq:diskbounds} yields
\begin{equation*}
 \beta(t)\le1-\biggl(\frac7{16}+\frac9{16}a^2\biggr)t+\frac9{16}a^2t^2\le1-\frac{31}{32}u+\frac9{16}u^2.
\end{equation*}
Since $0\le u\le1$ and $e^{-u/2}\ge1-\frac12 u+\frac18 u^2- \frac{1}{48} u^3$, we have
\begin{equation*}
 e^{-u/2}-\beta(t)\ge u\biggl(\frac{15}{32}-\frac7{16}u-\frac1{48}u^2\biggr)\ge\frac{u}{96}\ge0.
\end{equation*}

For $3\le s\le8$, we have $0<c_s\le1$ and $0\le c_su\le1$. Taylor's formula gives
\begin{equation*}
\begin{aligned}
 e^{-c_su}
 &\ge1-c_su+\frac{c_s^2u^2}{2}-\frac{c_s^3u^3}{6}+\frac{c_s^4u^4}{24}-\frac{c_s^5u^5}{120}\\
 &=1-c_su+\frac{11}{30}c_s^2u^2+\frac{c_s^2u^2}{120}(1-c_su)\bigl((c_su-2)^2+12\bigr)\ge1-c_su+\frac{11}{30}c_s^2u^2.
\end{aligned}
\end{equation*}
With $k=(s+1)/(2s)$ and $c_s$ from Table~\ref{tab:tail}, direct substitution gives $1-c_s+\frac{11}{30}c_s^2-k^2>0$ for $s=2,\,3,\,\ldots,\,8$. Since $0\le u\le1$, \eqref{eq:betau} implies
\begin{equation*}
\begin{aligned}
 e^{-c_su}-\beta(t)
 &\ge\biggl(1-c_su+\frac{11}{30}c_s^2u^2\biggr)
      -\bigl(1-u+k^2u^2\bigr)
 =u\biggl(1-c_s+\biggl(\frac{11}{30}c_s^2-k^2\biggr)u\biggr)\\
 &=u\biggl((1-c_s)(1-u)
       +\biggl(1-c_s+\frac{11}{30}c_s^2-k^2\biggr)u\biggr)
       \ge0,
\end{aligned}
\end{equation*}
which proves the lemma.
\end{proof}
\smallskip

\begin{proposition}\label{prop:tail}
Under the hypotheses of Theorem~\ref{thm:completeweighted}, suppose $2\le s\le8$ and $r\ge r_0(s)$.
Then $G$ has a zero $\zeta$ satisfying $\abs{a-\zeta}\le R_s$.
\end{proposition}
\begin{proof}
Suppose, for a contradiction, that $\abs{a-\zeta_j}>R_s$ for $j=1,\,2,\,\ldots,\,r$.
Substituting $u\=at$ in \eqref{eq:U} and applying Lemma~\ref{lem:taildecay}, we obtain
\begin{equation}\label{eq:etail}
\begin{aligned}
 E=rk^2\sum_{j=0}^{s-1}\frac{(rk)^j}{j!}
       \int_0^a\!u^{j+1}\beta\Bigl(\frac{u}{a}\Bigr)^{(r-1)/2}\,\mathrm{d}u
  &\le rk^2\sum_{j=0}^{s-1}\frac{(rk)^j}{j!}
       \int_0^{+\infty}\!u^{j+1}e^{-c_s(r-1)u/2}\,\mathrm{d}u\\
  &=\frac{4rk^2}{c_s^2(r-1)^2}
       \sum_{j=0}^{s-1}(j+1)
           \biggl(\frac{2rk}{c_s(r-1)}\biggr)^j.
\end{aligned}
\end{equation}
The last equality follows from $\int_0^{+\infty}\! u^{j+1}e^{-bu}\,\mathrm{d}u=(j+1)!/b^{j+2}$ for $b>0$.
The resulting upper bound decreases with $r$ for $r>1$, so it suffices to evaluate it at $r=r_0$.
Substituting the parameters in Table~\ref{tab:tail} into \eqref{eq:etail} gives the rational upper bounds in the following table.
\begin{table}[ht]
\centering\small
\begin{tabular}{c*{7}{c}}
\toprule
$s$ & 2 & 3 & 4 & 5 & 6 & 7 & 8\\
\midrule
upper bound for $E$ & $43/100$ & $39/100$ & $42/100$ & $54/100$
                   & $53/100$ & $53/100$ & $53/100$\\
\bottomrule
\end{tabular}
\par\vspace{8pt}
\small
\caption[Integral bounds]{The upper bounds are taken at $r=r_0(s)$, and the corresponding $c_s$ are taken as in Table~\ref{tab:tail}.}\label{tab:tailinitial}
\end{table}

For $M$, the bound $(s+k)(rk)^{s-1}k^r$ in \eqref{eq:Mbound} decreases for integer $r\ge r_0$ by \eqref{eq:degreeMonotonicity} with $j=s-1$, since $r_0\ge(s+1)^2+1\ge5s$. When $s=2$, we have $(s+k)(rk)^{s-1}k^r\le \frac{11}{4}\cdot\bigl(\frac{3}{4}\bigr)^{151}\cdot 150<\frac{1}{100}$. By a direct calculation, we also have $(s+k)(r_{0}k)^{s-1}k^{r_0}<\frac{1}{100}$ when $3\le s\le 8$. Thus, for every integer $r\ge r_0$, we have $U=M+E< \frac{1}{100}+\frac{59}{100}=\frac{3}{5}<1$, contradicting Lemma~\ref{lem:originineq}. Hence, the assertion follows.
\end{proof}

\subsubsection{Case IV}\label{sec:finiteanalytic}

Assume $2\le s\le8$ and $(s+1)^2<r<r_0(s)$. There are still 553 integer pairs $(s,r)$ left. For $s=2$, we divide $10\le r\le149$ into five consecutive ranges and choose the decay rate $\lambda$ separately for each range. For $3\le s\le8$, we use $\lambda=13/12$ for all remaining pairs. Under the counterexample hypothesis, we will show the strict bounds for $E$ and $M$ listed in the last two columns. In each row, the sum of the two upper bounds is at most $1$, so these estimates imply $U=M+E<1$, contradicting Lemma~\ref{lem:originineq}.

\begin{table}[ht]
\centering
\renewcommand{\arraystretch}{2.2}
\setlength{\tabcolsep}{6pt}
\setlength{\aboverulesep}{3pt}
\setlength{\belowrulesep}{3pt}
\begin{tabular}{c|cccc|cc}
\toprule
Row&$s$&$r_-$&$r_+$&$\lambda$&
$E<$&$M<$\\
\midrule
1&2&10&10&$\dfrac{13}{10}$&$\dfrac{59}{100}$&$\dfrac{37}{100}$\\
2&2&11&12&$\dfrac{5}{4}$&$\dfrac{58}{100}$&$\dfrac{37}{100}$\\
3&2&13&24&$\dfrac{11}{10}$&$\dfrac{67}{100}$&$\dfrac{30}{100}$\\
4&2&25&29&$1$&$\dfrac{41}{100}$&$\dfrac{3}{100}$\\
5&2&30&149&$\dfrac{7}{10}$&$\dfrac{90}{100}$&$\dfrac{1}{100}$\\
6&$3\le s\le8$&$(s+1)^2+1$&$r_0(s)-1$&$\dfrac{13}{12}$&
$\dfrac{88}{100}$&$\dfrac{12}{100}$\\
\bottomrule
\end{tabular}
\par\vspace{8pt}
\small
\caption[Finite-degree bounds]{$r_-$ and $r_+$ are the first and last integer degrees, respectively, and $\lambda$ is the
exponential decay rate which will be used in Lemma~\ref{lem:finiteexp}. The strict upper bounds in the last two columns will be obtained in Proposition~\ref{prop:finiteanalytic}.}
\label{tab:finiteanalytic}
\end{table}
\FloatBarrier

\begin{lemma}\label{lem:finiteexp}
Let $s,\,r,\,\lambda$ be as in Table~\ref{tab:finiteanalytic} and suppose $\abs{a-\zeta_j}>R_s$ for $j=1,\,2,\,\ldots,\,r$.
Then $\beta(t)\le e^{-\lambda at}$ for $0\le t\le1$.
\end{lemma}

\begin{proof}
Fix a row of Table~\ref{tab:finiteanalytic} and $s,\,r$ in its indicated range, taking $\lambda$ from that row.
Recall that $a$ is the distinguished point in $G(z)=(z-a)g'(z)+sg(z)$ and, by Lemma~\ref{lem:reduction}, $(s-1)/(s+1)<a<1$.
In the notation of \eqref{eq:notation} and \eqref{eq:betadefinition}, $k=(s+1)/(2s)\le 3/4$ and $\beta(t)=1-2axt+a^2k^2t^2$.
Since $x\le k$ by Lemma~\ref{lem:disk} and $ak<1$, we have $\beta(t)\ge(1-akt)^2>0$ for $0\le t\le1$.

By Lemma~\ref{lem:mean}, for $0\le t\le1$,
\begin{equation}\label{eq:finiteB}
 0<\beta(t)\le 1-\biggl(1+\frac{2}{2s+r(1-a^2)}-\bigl(1-a^2\bigr)k^2\biggr)t+a^2k^2t^2=:B_r(a, t).
\end{equation}

Fix $a$ and $r$.
We claim that $B_r(a,t)\le e^{-\lambda at}$ for $0\le t\le1$.
Then the lemma follows immediately from this claim and \eqref{eq:finiteB}.

Indeed, define $f(t)\=\log B_r(a,t)+\lambda at$ for $t\in[0, 1]$.
Then $f(0)=0$.
Since $r\le r_+$ and $1-a^2>0$, direct calculation gives
\begin{equation}\label{eq:finiteendconditions}
\begin{aligned}
 f'(0)&=\lambda a-1-\frac{2}{2s+r(1-a^2)}+\bigl(1-a^2\bigr)k^2\\
 &\le\lambda a-1-\frac{2}{2s+r_+(1-a^2)}+\bigl(1-a^2\bigr)k^2
 =-\frac{P_{s,r_+,\lambda}(a)}{2s+r_+(1-a^2)},
\end{aligned}
\end{equation}
where $P_{s,r_+,\lambda}(a)\= \bigl(1-\bigl(1-a^2\bigr)k^2-\lambda a\bigr)\bigl(2s+r_+\bigl(1-a^2\bigr)\bigr)+2$.

Recall that $S_6(-y)=\sum_{j=0}^{5}(-y)^j/j!$.
Since $5$ is odd, Taylor's formula with integral remainder gives $S_6(-y)<e^{-y}$ for $y>0$.
Hence
\begin{equation}\label{eq:finitefendpoint}
\begin{aligned}
 f(1)&=\log\biggl(k^2-\frac{2}{2s+r(1-a^2)}\biggr)+\lambda a
 \le\log\biggl(k^2-\frac{2}{2s+r_+(1-a^2)}\biggr)+\lambda a\\
 &=\log\biggl(S_6(-\lambda a)-\frac{Q_{s,r_+,\lambda}(a)}{2s+r_+(1-a^2)}\biggr)+\lambda a,
\end{aligned}
\end{equation}
where $Q_{s,r_+,\lambda}(a)\= \bigl(S_6(-\lambda a)-k^2\bigr)\bigl(2s+r_+\bigl(1-a^2\bigr)\bigr)+2$.

Appendix~\ref{app:finiteidentities} proves $P_{s,r_+,\lambda}(a)>0$ and $Q_{s,r_+,\lambda}(a)>0$ for $a\in[(s-1)/(s+1),1]$ by positive Bernstein expansions.
Thus $f'(0)<0$ by \eqref{eq:finiteendconditions}, and $S_6(-\lambda a)>k^2-\frac{2}{2s+r_+(1-a^2)} \ge k^2-\frac{2}{2s+r(1-a^2)}=B_r(a,1)>0$, so \eqref{eq:finitefendpoint} and the Taylor bound give $f(1)<\log S_6(-\lambda a)+\lambda a<0$.

Furthermore,
\begin{equation*}
\begin{aligned}
 B_r(a,t)f'(t)
 &=\lambda a-1-\frac{2}{2s+r(1-a^2)}+\bigl(1-a^2\bigr)k^2\\
 &\quad+\biggl(2a^2k^2-\lambda a
       \biggl(1+\frac{2}{2s+r(1-a^2)}-\bigl(1-a^2\bigr)k^2\biggr)\biggr)t
       +\lambda a^3k^2t^2.
\end{aligned}
\end{equation*}
Since the right-hand side has positive leading coefficient, $B_r(a,t)>0$, and $f'(0)<0$, $f$ is decreasing, or first decreasing and then increasing, on $[0,1]$.
Hence $f(t)\le\max\{f(0),\,f(1)\}=0$ on $[0,1]$, which proves the claim.
\end{proof}
\smallskip

\begin{proposition}\label{prop:finiteanalytic}
Under the hypotheses of Theorem~\ref{thm:completeweighted}, suppose $2\le s\le8$ and $(s+1)^2<r<r_0(s)$.
Then $G$ has a zero $\zeta$ satisfying $\abs{a-\zeta}\le R_s$.
\end{proposition}

\begin{proof}
Suppose, for a contradiction, that $\abs{a-\zeta_j}>R_s$ for $j=1,\,2,\,\ldots,\,r$.
Fix the row of Table~\ref{tab:finiteanalytic} containing $(s,r)$.
Substituting $u=at$ in \eqref{eq:U} and applying Lemma~\ref{lem:finiteexp}, we have
\begin{equation}\label{eq:finitee}
\begin{aligned}
 E&=rk^2\sum_{j=0}^{s-1}\frac{(rk)^j}{j!}
       \int_0^a\!u^{j+1}\beta\Bigl(\frac{u}{a}\Bigr)^{(r-1)/2}\,\mathrm{d}u\\
  &\le rk^2\sum_{j=0}^{s-1}\frac{(rk)^j}{j!}
       \int_0^{+\infty}\!u^{j+1}e^{-\lambda(r-1)u/2}\,\mathrm{d}u\\
  &=\frac{4k^2}{\lambda^2}
       \sum_{j=0}^{s-1}(j+1)
       \biggl(\frac{2k}{\lambda}\biggr)^j
       \frac{r^{j+1}}{(r-1)^{j+2}}\\
  &\le\frac{4k^2}{\lambda^2}
       \sum_{j=0}^{s-1}(j+1)
       \biggl(\frac{2k}{\lambda}\biggr)^j
       \frac{r_-^{j+1}}{(r_--1)^{j+2}}.
\end{aligned}
\end{equation}
The last inequality holds because $\frac{\mathrm d}{\mathrm dr}\bigl(\frac{r^{j+1}}{(r-1)^{j+2}}\bigr)<0$ for $r>1$ and $j=0,\,1,\,\ldots,\,s-1$.
Substituting the values of $s$, $\lambda$, and $r_-$ from Table~\ref{tab:finiteanalytic} into \eqref{eq:finitee}, with $k=(s+1)/(2s)$, yields the strict bounds for $E$ in Table~\ref{tab:finiteanalytic}.
The corresponding calculations are recorded in Table~\ref{tab:integralarithmetic} of the appendix.

For $M$, due to $a<1$ and $\beta(1)\le k^2$ in \eqref{eq:U}, we obtain that
\begin{equation}\label{eq:finitem}
 M\le S_s(rk)k^r+\frac{(rk)^s k^r}{(s-1)!(r+s)}
 =\sum_{j=0}^{s-1}\frac{k^j}{j!}\,r^jk^r
  +\frac{k^s}{(s-1)!}\frac{r^sk^r}{r+s}.
\end{equation}
Since $r\ge r_-\ge(s+1)^2+1\ge5s$, \eqref{eq:degreeMonotonicity} gives $r^jk^r\le r_-^jk^{r_-}$ for $j=0,\,1,\,\ldots,\,s$.
Together with $r+s\ge r_-+s$, we have
\begin{equation}\label{eq:finitemendpoint}
 M\le\sum_{j=0}^{s-1}\frac{k^j}{j!}\,r_-^jk^{r_-}
       +\frac{k^s}{(s-1)!}\frac{r_-^sk^{r_-}}{r_-+s}.
\end{equation}

Substituting $s$ and $r_-$ into this inequality gives the bounds for $M$ in rows 4--6 of Table~\ref{tab:finiteanalytic}.
The corresponding calculations are recorded in Table~\ref{tab:endpointarithmetic} of the appendix.

For rows 1--3, $s=2$, $a>(s-1)/(s+1)=1/3$, and $k=3/4$.
The sharper estimate for $\beta(1)$ in \eqref{eq:finiteB} gives
\begin{equation*}
 0<\beta(1)\le\frac9{16}-\frac{2}{4+r(1-a^2)}
 \le\frac9{16}-\frac{2}{4+8r_+/9}<1.
\end{equation*}
Then \eqref{eq:U} implies that
\begin{equation*}
 M\le\biggl(1+\frac{3r}{4}\biggr)
       \biggl(\frac9{16}-\frac{2}{4+8r_+/9}\biggr)^{\lfloor r/2\rfloor}
       +\frac{(3r/4)^2(3/4)^r}{r+2}.
\end{equation*}
Since $r\ge 10$,
\begin{equation*}
\begin{aligned}
 \frac{\bigl(1+\frac{3(r+2)}4\bigr)
       \bigl(\frac9{16}-\frac{2}{4+8r_+/9}\bigr)^{\lfloor(r+2)/2\rfloor}}
      {\bigl(1+\frac{3r}4\bigr)
       \bigl(\frac9{16}-\frac{2}{4+8r_+/9}\bigr)^{\lfloor r/2\rfloor}}
 &=\biggl(1+\frac{3/2}{1+3r/4}\biggr)
       \biggl(\frac9{16}-\frac{2}{4+8r_+/9}\biggr)\\
 &\le\frac9{16}\biggl(1+\frac{3/2}{1+3r/4}\biggr)
 \le\frac9{16}\frac{20}{17}<1  \quad\text{ and} \\
 \frac{(3(r+1)/4)^2(3/4)^{r+1}(r+2)}{(3r/4)^2(3/4)^r(r+3)}
 &=\frac34\biggl(1+\frac1r\biggr)^2\frac{r+2}{r+3}
 \le\frac34\biggl(\frac{11}{10}\biggr)^2
 =\frac{363}{400}<1.
\end{aligned}
\end{equation*}
Consequently,
\begin{equation}\label{eq:finitemsharp}
 M\le
 \biggl(1+\frac{3(r_-+1)}4\biggr)
 \biggl(\frac9{16}-\frac{2}{4+8r_+/9}\biggr)^{\lfloor r_-/2\rfloor}
 +\frac{(3r_-/4)^2(3/4)^{r_-}}{r_-+2}.
\end{equation}
Substituting $r_+$ and $r_-$ into \eqref{eq:finitemsharp} gives the bounds for $M$ in rows 1--3 of Table~\ref{tab:finiteanalytic}.
The corresponding calculations are recorded in Table~\ref{tab:sharpendpointarithmetic} of the appendix.

Thus $U=M+E<1$, contradicting Lemma~\ref{lem:originineq}. Therefore, the assertion follows.
\end{proof}

\smallskip

\subsection{Proof of Theorem~\ref{thm:completeweighted}}\label{sec:generalassembly}
\begin{proof}
Propositions~\ref{prop:small}, \ref{prop:large}, \ref{prop:tail}, and \ref{prop:finiteanalytic} give the conclusion in all of the four parameter ranges listed in Table~\ref{tab:cases}, which completes the proof.
\end{proof}
\smallskip

\section{Proof of Theorem~\ref{prop:fourpoints}}\label{sec:quartic}

Let $\lambda,\,w_2,\,w_3,\,w_4$ satisfy the hypotheses of Theorem~\ref{prop:fourpoints}. Write $\mathcal B(0, w_2, w_3, w_4)=\overline{D}(c_w,R)$, where $c_w$ is the center and $R$ is the radius. We may assume $R^2\le3/4$, since otherwise the required lower bound is immediate.

\subsection{Centered coordinates}

For $j=2,\,3,\,4$, denote
\begin{equation}\label{eq:quartic-center}
 p_t\=\frac{2(1+\lambda)}3,\qquad
 t_j\=p_t-w_j,\qquad c_t\=p_t-c_w,\qquad
 u\=(2\lambda-1)(\lambda-2).
\end{equation}
Then $\mathcal B(p_t,\,t_2,\,t_3,\,t_4) =\overline{D}(c_t,R)$.

\begin{lemma}\label{lem:quartic-centering}
Suppose $R^2\le3/4$. Then $\sum_{j=2}^4t_j=0$ and $\sum_{j=2}^4t_j^2=4u/3$. Moreover, $\abs{\lambda}<3/2$ and $0<\abs{u}<27/16$.
\end{lemma}
\begin{proof}
By~\eqref{eq:quartic-sums}, $\sum_{j=2}^4w_j=3p_t$, and hence $\sum_{j=2}^4t_j=0$. Furthermore,
\begin{equation*}
\begin{aligned}
 \sum_{j=2}^4t_j^2
 &=\sum_{j=2}^4w_j^2-2p_t\sum_{j=2}^4w_j+3p_t^2
 =\sum_{j=2}^4w_j^2-3p_t^2\\
 &=\Biggl(\sum_{j=2}^4w_j\Biggr)^2
   -2\sum_{2\le i<j\le4}w_iw_j-3p_t^2
 =6p_t^2-12\lambda
  =\frac43(2\lambda-1)(\lambda-2),
\end{aligned}
\end{equation*}
so $\sum_{j=2}^4t_j^2=4u/3$, and hence
\begin{equation}\label{eq:quartic-mean}
 3R^2\ge\sum_{j=2}^4\abs{t_j-c_t}^2
 =\sum_{j=2}^4\abs{t_j}^2+3\abs{c_t}^2
 \ge\sum_{j=2}^4\abs{t_j}^2
 \ge\AbsBigg{\sum_{j=2}^4t_j^2}
 =\frac43\abs{u},
\end{equation}
which gives $R^2\ge4\abs{u}/9$.

Since $\abs{c_w}\le R$ and $\abs{w_j}\le2R$,
\begin{equation*}
 6\abs{\lambda}
 =\Bigl\lvert\sum_{2\le i<j\le4}w_iw_j\Bigr\rvert
 \le\sum_{2\le i<j\le4}\abs{w_iw_j}\le12R^2\le9.
\end{equation*}
If $\abs{\lambda}=3/2$, then $R^2=3/4$, $\abs{w_j}=2R$, and hence
\begin{equation*}
 2R=\abs{w_j}\le\abs{w_j-c_w}+\abs{c_w}\le2R,
\end{equation*}
for $j=2,\,3,\,4$.
Thus $w_2=w_3=w_4=2c_w$.
Then $t_2=t_3=t_4=0$ by \eqref{eq:quartic-center} and $\sum_{j=2}^4t_j=0$.
It follows from $\sum_{j=2}^4t_j^2=4u/3$ that $u=0$, so $\lambda\in\{1/2,\,2\}$, contradicting $\abs{\lambda}=3/2$.
Since $u=(2\lambda-1)(\lambda-2)$ and $1\le\abs{\lambda}<3/2$, we have $u\ne0$.
If $\abs{u}=27/16$, equality in~\eqref{eq:quartic-mean} gives $c_t=0$, $\abs{t_j}=R$, and a common argument for the three $t_j^2$.
Thus, $t_j=\varepsilon_jR e^{\I\theta}$ with $\varepsilon_j\in\{-1,\,1\}$, contradicting $\sum_{j=2}^4t_j=0$.

Therefore, the conclusion holds.
\end{proof}

\smallskip

\subsection{Ellipse representation}
By Lemma~\ref{lem:quartic-centering}, $u\ne0$ and $\sum_{j=2}^4t_j^2=4u/3$.
Choose $\eta\in\R$ so that $e^{-2\I\eta}u=-\abs{u}$.
For $j=2,\,3,\,4$, denote
\begin{equation}\label{eq:quartic-rotation}
 p_\tau\=e^{-\I\eta}p_t\=\alpha+\I\beta,\qquad
 \tau_j\=e^{-\I\eta}t_j,\qquad
 c_\tau\=e^{-\I\eta}c_t\=\gamma+\I\delta,
\end{equation}
where $\alpha,\,\beta,\,\gamma,\,\delta\in\R$. Then $\mathcal B(p_\tau, \tau_2, \tau_3, \tau_4)=\overline{D}(c_\tau,R)$. Figure~\ref{fig:reflection} below illustrates the transformations of the coordinates.

Set $\nu\=e^{2\pi \I/3}$.
For $z_2,\,z_3,\,z_4\in\C$ with $\sum_{j=2}^4z_j=0$, the discrete Fourier transform formulas with $n=3$ give $z_j=\Pi_+\nu^j+\Pi_-\nu^{-j}$ for $j=2,\,3,\,4$, where $\Pi_\pm\=\frac13\sum_{k=2}^4z_k\nu^{\mp k}$~\cite[(3.11.38)--(3.11.39)]{NISTDLMF}.
Since $1+\nu+\nu^2=0$, this also implies $\sum_{j=2}^4z_j^2=6\Pi_+\Pi_-$.

\begin{lemma}\label{lem:quartic-ellipsebounds}
After relabeling, there exist $a_0,\,b_0,\,\psi\in\R$ satisfying $0\le a_0<b_0$, $b_0^2-a_0^2=8\abs{u}/9$, and
\begin{equation}\label{eq:quartic-ellipse}
 \tau_j=a_0\cos(\psi+(2/3)\pi j)+\I b_0\sin(\psi+(2/3)\pi j),
\end{equation}
for $j=2,\,3,\,4$. Let $K\=\sqrt{\abs{c_\tau}^2a_0^2+8\abs{u}\delta^2/9}$. Then
\begin{equation}\label{eq:quartic-radius}
 R^2\ge\max\bigl\{
 (\abs{\gamma}-\abs{\alpha})^2
 +(\abs{\delta}-\abs{\beta})^2,\,
 a_0^2+ (2/9)\abs{u}+\abs{c_\tau}^2+K,\,
 a_0^2+ (2/3)\abs{u}+\abs{c_\tau}^2-K
 \bigr\}.
\end{equation}
\end{lemma}

\begin{remark}
Recall that
\begin{equation*}
 R=\min_{c\in\C}\max\Bigl\{\abs{c},\,
       \max_{2\le j\le4}\abs{w_j-c}\Bigr\}
  =\max\Bigl\{\abs{c_w},\,
       \max_{2\le j\le4}\abs{w_j-c_w}\Bigr\}.
\end{equation*}
The first term in~\eqref{eq:quartic-radius} comes from requiring the disk $\overline{D}(c_w,R)$ to contain $0$, namely $R^2\ge\abs{c_w}^2=\abs{p_\tau-c_\tau}^2$. The last two come from requiring the same disk to contain all three $w_j$, namely $R^2\ge\max_{2\le j\le4}\abs{w_j-c_w}^2 =\max_{2\le j\le4}\abs{\tau_j-c_\tau}^2$.
\end{remark}
\begin{proof}[Proof of Lemma~\ref{lem:quartic-ellipsebounds}]
Since $\sum_{j=2}^4\tau_j=0$, the inversion formula with $z_j=\tau_j$ and Lemma~\ref{lem:quartic-centering} yield
\begin{equation*}
 \Pi_+\Pi_-=\frac16\sum_{j=2}^4\tau_j^2
 =\frac{e^{-2\I\eta}}6\sum_{j=2}^4t_j^2
 =\frac29e^{-2\I\eta}u=-\frac{2\abs{u}}9.
\end{equation*}
Interchange $\tau_2$ and $\tau_4$ if necessary so that $\abs{\Pi_+}\ge\abs{\Pi_-}$, and define $a_0\=\abs{\Pi_+}-\abs{\Pi_-}$ and $b_0\=\abs{\Pi_+}+\abs{\Pi_-}$. Then $0\le a_0<b_0$ and $b_0^2-a_0^2=8\abs{u}/9$. Since $\Pi_+\Pi_-<0$, there exists $\psi\in\R$ satisfying
\begin{equation*}
 \Pi_+=\frac{a_0+b_0}{2}e^{\I\psi}\quad\text{ and }\quad
 \Pi_-=\frac{a_0-b_0}{2}e^{-\I\psi}.
\end{equation*}
Substituting these expressions into the discrete Fourier transform formula $\tau_j=\Pi_+\nu^j+\Pi_-\nu^{-j}$ yields~\eqref{eq:quartic-ellipse}.

The first lower bound in~\eqref{eq:quartic-radius} follows from $R^2\ge\abs{p_\tau-c_\tau}^2$. For the other two lower bounds, note that $\Rea\bigl(e^{2\I\psi}\nu^{2j} \bigr)=\Rea\bigl(e^{-2\I\psi}\nu^j\bigr)$ for $j=2,\,3,\,4$.
Then
\begin{equation*}
\begin{aligned}
 \abs{\tau_j-c_\tau}^2
 &=\abs{\tau_j}^2-2\Rea(\overline{c_\tau}\tau_j)+\abs{c_\tau}^2\\
 &=\frac{a_0^2+b_0^2}{2}
   +\frac{a_0^2-b_0^2}{2}\Rea\bigl(e^{2\I\psi}\nu^{2j}\bigr)
   -2\Rea\bigl((\gamma a_0-\I\delta b_0)e^{\I\psi}\nu^j\bigr)
   +\abs{c_\tau}^2\\
 &=\frac{a_0^2+b_0^2}{2}+\abs{c_\tau}^2
   +\Rea\biggl(\biggl(\frac{a_0^2-b_0^2}{2}e^{-2\I\psi}
   -2(\gamma a_0-\I\delta b_0)e^{\I\psi}\biggr)\nu^j\biggr)\\
 &=a_0^2+\frac{4\abs{u}}9+\abs{c_\tau}^2-\Rea\bigl(\Gamma\nu^j\bigr),
\end{aligned}
\end{equation*}
where $\Gamma\=\frac{4\abs{u}}9e^{-2\I\psi} +2(\gamma a_0-\I\delta b_0)e^{\I\psi}$.
For $\Gamma\ne0$, one of the three $\Gamma\nu^j$ has argument within $\pi/3$ of $\pi$, so $\min_{j\in\{2,\,3,\,4\}}\Rea\bigl(\Gamma\nu^j\bigr)\le-\abs{\Gamma}/2$.
The same inequality holds for $\Gamma=0$.
Moreover, $\abs{\gamma a_0-\I\delta b_0}^2 =\abs{c_\tau}^2a_0^2+8\abs{u}\delta^2/9=K^2$, so $\abs{\Gamma}\ge2\abs{K-2\abs{u}/9}$.
Consequently,
\begin{equation*}
 R^2\ge a_0^2+\frac{4\abs{u}}9+\abs{c_\tau}^2+\frac12\abs{\Gamma}
 \ge a_0^2+\frac{4\abs{u}}9+\abs{c_\tau}^2 +\Absbigg{K-\frac{2\abs{u}}{9}},
\end{equation*}
which yields the remaining two lower bounds in~\eqref{eq:quartic-radius}.
\end{proof}

\smallskip

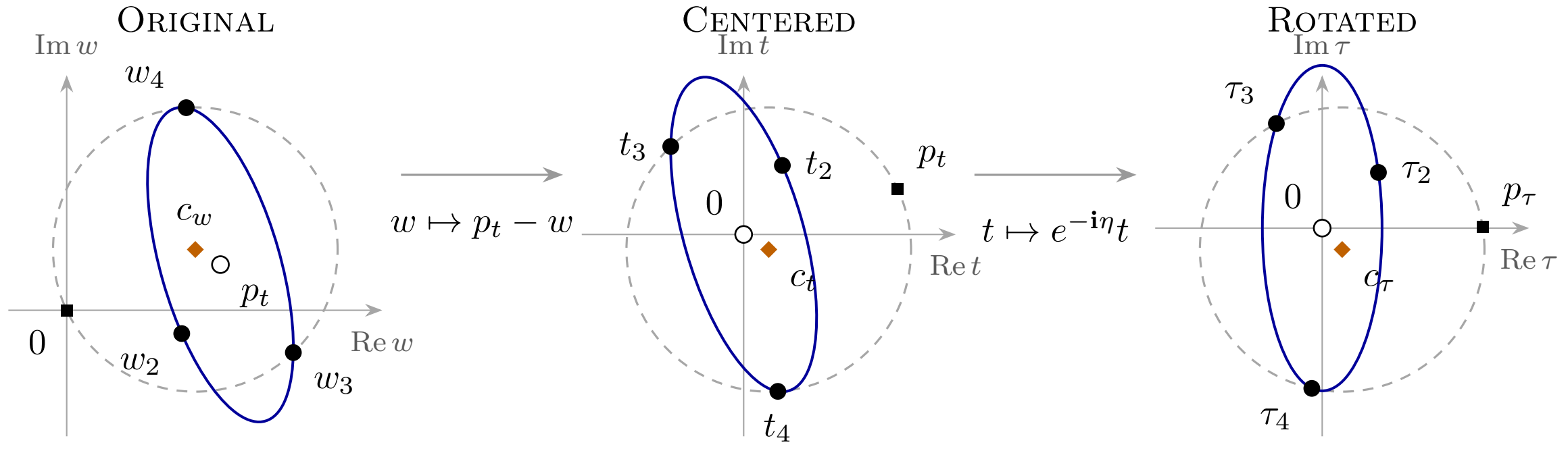
\begin{figure}[!htbp]
\centering
\begin{tikzpicture}[scale=1.22,every node/.style={font=\small},>=Stealth,
 ellipse curve/.style={blue!60!black,line width=.75pt},
 axis/.style={->,black!35,line width=.4pt},
 axislabel/.style={font=\scriptsize,text=black!65},
 disk curve/.style={gray!70,dash pattern=on 3pt off 2.5pt,line width=.6pt},
 root/.style={circle,fill=black,inner sep=1.65pt},
 distinguished/.style={rectangle,fill=black,inner sep=1.65pt},
 centroid/.style={circle,draw=black,fill=white,inner sep=1.65pt,line width=.5pt}]
\begin{scope}
 \node[font=\small\scshape] at (0,1.85) {Original};
 \draw[axis] (-1.5,-.487782)--(1.5,-.487782)
   node[axislabel,below=2pt] {$\Rea w$};
 \draw[axis] (-1.031901,-1.5)--(-1.031901,1.4)
   node[axislabel,above=2pt] {$\operatorname{Im} w$};
 \draw[disk curve] (0,0) circle (1.141381);
 \draw[ellipse curve,rotate around={15.984421:(.201433,-.121115)}]
   (.201433,-.121115) ellipse (.48 and 1.306829);
 \node[distinguished,label={below left:$0$}] at (-1.031901,-.487782) {};
 \node[root,label={below left:$w_2$}] at (-.109098,-.675005) {};
 \node[root,label={below right:$w_3$}] at (.786237,-.827395) {};
 \node[root,label={above left:$w_4$}] at (-.072841,1.139054) {};
 \node[centroid,label={below right:$p_t$}] at (.201433,-.121115) {};
 \draw[orange!75!black,fill=orange!75!black]
   (0,.06)--(.06,0)--(0,-.06)--(-.06,0)--cycle;
 \node[above=3pt] at (0,0) {$c_w$};
\end{scope}
\draw[->,gray!80,line width=.6pt] (1.65,.6)--(2.95,.6);
\node[below] at (2.3,.42) {$w\mapsto p_t-w$};
\begin{scope}[xshift=4.6cm]
 \node[font=\small\scshape] at (0,1.85) {Centered};
 \draw[axis] (-1.5,.121115)--(1.5,.121115)
   node[axislabel,below=2pt] {$\Rea t$};
 \draw[axis] (-.201433,-1.5)--(-.201433,1.4)
   node[axislabel,above=2pt] {$\operatorname{Im} t$};
 \draw[disk curve] (0,0) circle (1.141381);
 \draw[ellipse curve,rotate around={15.984421:(-.201433,.121115)}]
   (-.201433,.121115) ellipse (.48 and 1.306829);
 \node[distinguished,label={above right:$p_t$}] at (1.031901,.487782) {};
 \node[root,label={right:$t_2$}] at (.109098,.675005) {};
 \node[root,label={left:$t_3$}] at (-.786237,.827395) {};
 \node[root,label={below:$t_4$}] at (.072841,-1.139054) {};
 \node[centroid,label={above left:$0$}] at (-.201433,.121115) {};
 \draw[orange!75!black,fill=orange!75!black]
   (0,.06)--(.06,0)--(0,-.06)--(-.06,0)--cycle;
 \node[below right=2pt] at (0,0) {$c_t$};
\end{scope}
\draw[->,gray!80,line width=.6pt] (6.25,.6)--(7.55,.6);
\node[below] at (6.9,.42) {$t\mapsto e^{-\I\eta}t$};
\begin{scope}[xshift=9.2cm]
 \node[font=\small\scshape] at (0,1.85) {Rotated};
 \draw[axis] (-1.5,.171902)--(1.5,.171902)
   node[axislabel,below=2pt] {$\Rea\tau$};
 \draw[axis] (-.160292,-1.5)--(-.160292,1.4)
   node[axislabel,above=2pt] {$\operatorname{Im}\tau$};
 \draw[disk curve] (0,0) circle (1.141381);
 \draw[ellipse curve] (-.160292,.171902) ellipse (.48 and 1.306829);
 \node[distinguished,label={above right:$p_\tau$}] at (1.126327,.184762) {};
 \node[root,label={right:$\tau_2$}] at (.290760,.618864) {};
 \node[root,label={above left:$\tau_3$}] at (-.527994,1.011916) {};
 \node[root,label={below left:$\tau_4$}] at (-.243643,-1.115073) {};
 \node[centroid,label={above left:$0$}] at (-.160292,.171902) {};
 \draw[orange!75!black,fill=orange!75!black]
   (0,.06)--(.06,0)--(0,-.06)--(-.06,0)--cycle;
 \node[below right=2pt] at (0,0) {$c_\tau$};
\end{scope}
\end{tikzpicture}
\caption{Centering and rotation of the four-point configuration (schematic).
Each dashed circle has radius $R$.}
\label{fig:reflection}
\end{figure}
\FloatBarrier

\subsection{Some algebraic inequalities}
Define
\begin{equation*}
 J\=\frac92\abs{p_\tau}^2-2\abs{u}
       -\frac32\sqrt{2\abs{u}}\abs{\beta}\quad\text{ and }\quad
 \mathcal M\=\frac{\abs{u}}2
 \bigl(4J^2-81\alpha^2(3-2\abs{u})\bigr).
\end{equation*}
Consider the algebraic conditions
\begin{subequations}\label{eq:quartic-conditions}
\begin{align}
 J&>0,\label{eq:quartic-Jpositive}\\
 2J\biggl(1-\frac{3\abs{\beta}}{\sqrt{2\abs{u}}}\biggr)&<9\alpha^2,
   \label{eq:quartic-weight2}\\
 2J\biggl(1+\frac{3\abs{\beta}}{\sqrt{2\abs{u}}}\biggr)&<27\alpha^2,
   \label{eq:quartic-weight3}\\
 \abs{p_\tau}^2&>\frac34,\label{eq:quartic-pnorm}\\
 \biggl(\abs{p_\tau}-\frac{\sqrt3}2\biggr)^2
 &>\frac32-\frac43\abs{u},\label{eq:quartic-pdistance}\\
 \mathcal M&\ge0.\label{eq:quartic-Mnonnegative}
\end{align}
\end{subequations}
These conditions will be used in Lemma~\ref{lem:quartic-elimination}. Conditions~\eqref{eq:quartic-Jpositive}--\eqref{eq:quartic-weight3} ensure that the weights in~\eqref{eq:quartic-optimal-weights} are positive. Under $R^2\le3/4$, conditions~\eqref{eq:quartic-pnorm} and \eqref{eq:quartic-pdistance} allow us to simplify the radius bounds in~\eqref{eq:quartic-radius}. Finally, condition~\eqref{eq:quartic-Mnonnegative} ensures that the resulting lower bound for $R^2$ is at least $3/4$.

\begin{lemma}\label{lem:quartic-algebra}
Let $\lambda\in\C$ satisfy $1\le\abs{\lambda}<3/2$ and $0<\abs{u}<27/16$. Then the algebraic conditions~\eqref{eq:quartic-conditions} hold. Moreover, $\mathcal M=0$ if and only if $\lambda=1$.
\end{lemma}
\begin{proof}
By applying the Joukowski parametrization, we have
\begin{equation}\label{eq:quartic-joukowski}
 \lambda=\frac{(3\kappa+1)(\kappa+3)}{8\kappa}
 =\frac54+\frac38\biggl(\kappa+\frac{1}{\kappa}\biggr).
\end{equation}
The corresponding quadratic has roots of product one, so choose $\abs{\kappa}\ge1$ and write $\kappa=x+\I y$ with $x,\,y\in\R$. By~\eqref{eq:quartic-joukowski}, $1+\lambda=3\bigl(\kappa^2+6\kappa+1\bigr)\big/(8\kappa)$ and $u=(2\lambda-1)(\lambda-2)=9\bigl(\kappa^2-1\bigr)^2\big/\bigl(32\kappa^2\bigr)$. Since $u\ne0$, we have $\kappa\ne\pm1$. Set
\begin{equation*}
 H_1\=\Absbig{\kappa^2+6\kappa+1}^2\quad\text{ and }\quad
 H_2\=\Absbig{\kappa^2-1}^2>0.
\end{equation*}
Let $\alpha_0,\,\beta_0\in\R$ be given by $\alpha_0+\I\beta_0=\bigl(\kappa^2+6\kappa+1 \bigr)\overline{(\kappa^2-1)}$. Then
\begin{equation*}
 \alpha_0=\bigl(\abs{\kappa}^2-1\bigr)\bigl(\abs{\kappa}^2+6x+1\bigr)\quad\text{ and }\quad
 \beta_0=-2y\bigl(3\abs{\kappa}^2+2x+3\bigr).
\end{equation*}
By~\eqref{eq:quartic-joukowski}, $0\le64\abs{\kappa}^2\bigl(\abs{\lambda}^2-1\bigr) =3\bigl(\abs{\kappa}^2+6x+1\bigr) \bigl(3\abs{\kappa}^2+2x+3\bigr)$. Since $3\abs{\kappa}^2+2x+3=3\bigl(x+\frac13\bigr)^2+3y^2+\frac83>0$, we have $\abs{\kappa}^2+6x+1\ge0$, and hence $\alpha_0\ge0$. Moreover, since $e^{-2\I\eta}=-\overline u/\abs{u}$, we have $p_\tau^2=-\frac{4(1+\lambda)^2\overline u}{9\abs{u}} =-\frac{(\alpha_0+\I\beta_0)^2}{16\abs{\kappa}^2H_2}$. Then
\begin{equation}\label{eq:quartic-joukowski-invariants}
 \abs{p_\tau}^2=\frac{H_1}{16\abs{\kappa}^2},\qquad
 \abs{u}=\frac{9H_2}{32\abs{\kappa}^2},\qquad
 \alpha^2=\frac{\beta_0^2}{16\abs{\kappa}^2H_2},\qquad
 \beta^2=\frac{\alpha_0^2}{16\abs{\kappa}^2H_2}.
\end{equation}

To prove~\eqref{eq:quartic-Jpositive}--\eqref{eq:quartic-weight3}, let $\mathcal N_1$, $\mathcal N_2$, and $\mathcal N_3$ be the real polynomials in $x,\,y$ given by
\[\mathcal N_1\=H_1-2H_2-\alpha_0, \quad
 \mathcal N_2\=\beta_0^2-\mathcal N_1(H_2-\alpha_0), \quad\text{ and }\quad
 \mathcal N_3\=3\beta_0^2-\mathcal N_1(H_2+\alpha_0).\]
Here $x^2+y^2\ge1$, $1\le\abs{\lambda}<3/2$, and $0<\abs{u}<27/16$, with $\lambda$ given by~\eqref{eq:quartic-joukowski}.
On this parameter region, the identities~\eqref{eq:quartic-C1}--\eqref{eq:quartic-C3} in Appendix~\ref{app:quartic} give $\mathcal N_1,\,\mathcal N_2,\,\mathcal N_3>0$.

By~\eqref{eq:quartic-joukowski-invariants} and $\alpha_0\ge0$, $\sqrt{2\abs{u}}\abs{\beta}=3\alpha_0 \big/\bigl(16\abs{\kappa}^2\bigr)$.
Substituting~\eqref{eq:quartic-joukowski-invariants} and this identity into the definition of $J$ gives
\begin{equation*}
 J=\frac9{32\abs{\kappa}^2}(H_1-2H_2-\alpha_0)
 =\frac{9\mathcal N_1}{32\abs{\kappa}^2}>0,
\end{equation*}
which proves~\eqref{eq:quartic-Jpositive}. Since $3\abs{\beta}/\sqrt{2\abs{u}}=\alpha_0/H_2$, we have
\begin{align*}
9\alpha^2-2J\biggl(1-\frac{3\abs{\beta}}{\sqrt{2\abs{u}}}\biggr)
 &=\frac9{16\abs{\kappa}^2H_2}\mathcal N_2>0  \quad\text{ and}  \\
 27\alpha^2-2J\biggl(1+\frac{3\abs{\beta}}{\sqrt{2\abs{u}}}\biggr)
 &=\frac9{16\abs{\kappa}^2H_2}\mathcal N_3>0,
\end{align*}
which prove~\eqref{eq:quartic-weight2} and~\eqref{eq:quartic-weight3}.

Set $W\=H_1+6H_2-12\abs{\kappa}^2$. Identities~\eqref{eq:quartic-C6}--\eqref{eq:quartic-C5} give $H_1-12\abs{\kappa}^2>0$ and $W>4\sqrt{3H_1\abs{\kappa}^2}$, so
\begin{align*}
\abs{p_\tau}^2-\frac34
&=\frac1{16\abs{\kappa}^2}(H_1-12\abs{\kappa}^2)>0  \quad\text{ and}  \\
\biggl(\abs{p_\tau}-\frac{\sqrt3}2\biggr)^2 -\frac32+\frac43\abs{u}
&=\frac{W-4\sqrt{3H_1\abs{\kappa}^2}}{16\abs{\kappa}^2}>0,
\end{align*}
which prove~\eqref{eq:quartic-pnorm} and~\eqref{eq:quartic-pdistance}.

Finally,~\eqref{eq:quartic-Mcertificate} in Appendix~\ref{app:quartic} gives $\mathcal M\ge0$, which proves ~\eqref{eq:quartic-Mnonnegative}.

If $\mathcal M=0$, then~\eqref{eq:quartic-Mcertificate} gives $\abs{\kappa}^2=1$ and $x=-1/3$, so \eqref{eq:quartic-joukowski} gives $\lambda=\frac54+\frac34 x=1$. Conversely, $\lambda=1$ gives $\abs{u}=1$, $\alpha^2=16/9$, $\beta=0$, and $J=6$, so $\mathcal M=0$.
\end{proof}

\smallskip

\subsection{Radius estimates}\label{sec:quartic-elimination}
We now combine the bounds in Lemma~\ref{lem:quartic-ellipsebounds} to obtain a lower bound for $R^2$ independent of $a_0$ and $c_\tau$.
\begin{lemma}\label{lem:quartic-elimination}
Under the hypotheses of Theorem~\ref{prop:fourpoints}, $R^2\ge3/4$.
\end{lemma}
\begin{proof}
To prove the lemma, it suffices to consider $R^2\le3/4$ and show that $R^2=3/4$.
Lemma~\ref{lem:quartic-centering} gives $\abs{\lambda}<3/2$ and $0<\abs{u}<27/16$, so Lemma~\ref{lem:quartic-algebra} yields~\eqref{eq:quartic-conditions}. We now show that $\abs{c_\tau}^2\le\frac43\sqrt{2\abs{u}}\abs{\delta}$.
For $\xi\ge0$, define
\begin{equation*}
 \chi(\xi)\=\xi+ (2/3)\abs{u} + \abs{c_\tau}^2
 -\sqrt{\abs{c_\tau}^2\xi+ (8/9)\abs{u}\delta^2}.
\end{equation*}
The third term in~\eqref{eq:quartic-radius} gives $R^2\ge\chi\bigl(a_0^2\bigr)$.
When $c_\tau\ne0$ and $\xi>0$,
\begin{equation*}
 \chi'(\xi)=1-\frac{\abs{c_\tau}^2}
 {2\sqrt{\abs{c_\tau}^2\xi+(8/9)\abs{u}\delta^2}}\quad\text{ and }\quad
 \chi''(\xi)=\frac{\abs{c_\tau}^4}
 {4(\abs{c_\tau}^2\xi+(8/9)\abs{u}\delta^2)^{3/2}}>0.
\end{equation*}
For $c_\tau=0$, $\chi(\xi)=\xi+(2/3)\abs{u}$ is increasing.
Thus $\chi$ is convex on $[0,+\infty)$.

If $\abs{c_\tau}^2>\frac43\sqrt{2\abs{u}}\abs{\delta}$, then $\chi'(\xi)=0$ at $\xi=\frac{\abs{c_\tau}^2}{4}-\frac{8\abs{u}\delta^2}{9\abs{c_\tau}^2} >0$, and
\begin{equation*}
 R^2\ge\chi\bigl(a_0^2\bigr)\ge\min_{\xi\ge0}\chi(\xi)
 =\frac{2\abs{u}}3+\frac{3\abs{c_\tau}^2}4
  -\frac{8\abs{u}\delta^2}{9\abs{c_\tau}^2}
 >\frac{2\abs{u}}3+\frac{\abs{c_\tau}^2}2.
\end{equation*}
Hence $\abs{c_\tau}^2<\frac32-\frac43\abs{u}$.
Since~\eqref{eq:quartic-pnorm} gives $\abs{p_\tau}>\sqrt3/2$, the condition $\abs{p_\tau-c_\tau}\le R\le\sqrt3/2$ also yields $\abs{c_\tau}\ge\abs{p_\tau}-\frac12\sqrt3>0$.
Thus,
$ \bigl(\abs{p_\tau}-\frac{\sqrt3}2\bigr)^2
 \le\abs{c_\tau}^2<\frac32-\frac{4\abs{u}}3$.
This contradicts~\eqref{eq:quartic-pdistance}.
Hence $\abs{c_\tau}^2\le\frac43\sqrt{2\abs{u}}\abs{\delta}$.

If $\delta\ne0$, then $\chi'(0)=1-3\abs{c_\tau}^2 \bigl(4\sqrt{2\abs{u}}\abs{\delta}\bigr)^{-1}\ge0$, so the minimum is $\chi(0)$.
If $\delta=0$, then $c_\tau=0$ and the same conclusion holds.
Together with $a_0^2\ge0$ and $K\ge2\sqrt{2\abs{u}}\abs{\delta}/3$ in the second term of~\eqref{eq:quartic-radius}, this gives
\begin{equation}\label{eq:quartic-simplified}
\begin{aligned}
 R^2&\ge a_0^2+ (2/9)\abs{u}+\abs{c_\tau}^2+K
 \ge\abs{c_\tau}^2+(2/3)\sqrt{2\abs{u}} \abs{\delta}
      +(2/9)\abs{u} \quad\text{ and} \\
 R^2&\ge\chi\bigl(a_0^2\bigr)\ge\chi(0)
 =\abs{c_\tau}^2-(2/3)\sqrt{2\abs{u}} \abs{\delta}
      +(2/3)\abs{u} .
\end{aligned}
\end{equation}

For $\mu_1,\,\mu_2,\,\mu_3\ge0$ with $\mu_1+\mu_2+\mu_3=1$, define
\begin{equation}\label{eq:quartic-quadratic}
\begin{aligned}
 Q_\mu(X,Y)\={}&\mu_1\bigl((X-\abs{\alpha})^2
 +(Y-\abs{\beta})^2\bigr)\\
 &+\mu_2\bigl(X^2+Y^2+(2/3)\sqrt{2\abs{u}} Y
                      +(2/9)\abs{u} \bigr)\\
 &+\mu_3\bigl(X^2+Y^2-(2/3)\sqrt{2\abs{u}} Y
                      +(2/3)\abs{u} \bigr)\\
 ={}&(X-\mu_1\abs{\alpha})^2
 +\bigl(Y-\mu_1\abs{\beta} + \bigl(\sqrt{2\abs{u}}\big/3\bigr) (\mu_2-\mu_3)\bigr)^2
 +\Delta(\mu),
\end{aligned}
\end{equation}
where $X,\,Y\in\R$ and
\begin{equation*}
\Delta(\mu)
\=\mu_1\abs{p_\tau}^2+(2/9)\abs{u}\mu_2+(2/3)\abs{u}\mu_3
-\mu_1^2\alpha^2-\bigl(\mu_1\abs{\beta}- \bigl(\sqrt{2\abs{u}} \big/3\bigr) (\mu_2-\mu_3)\bigr)^2.
\end{equation*}
Combining the first term of~\eqref{eq:quartic-radius} with~\eqref{eq:quartic-simplified} yields
\begin{equation*}
 R^2\ge Q_\mu(\abs{\gamma},\abs{\delta})
 \ge\min \bigl\{Q_\mu(X,Y) : (X,Y)\in\R^2 \bigr\}
 =\Delta(\mu).
\end{equation*}
Now it suffices to choose the weights such that $\Delta(\mu)\ge3/4$.

By applying $\mu_2+\mu_3=1-\mu_1$ and completing the square in $\mu_2-\mu_3$, we obtain
\begin{equation*}
\begin{aligned}
 \Delta={}&\frac{\abs{u}}2
 +\mu_1\biggl(\abs{p_\tau}^2-\frac{4\abs{u}}9
              -\frac{\sqrt{2\abs{u}}}{3}\abs{\beta}\biggr)
 -\mu_1^2\alpha^2
 -\frac{2\abs{u}}9
 \biggl(\mu_2-\mu_3
       -\frac{3\mu_1\abs{\beta}}{\sqrt{2\abs{u}}}
       +\frac12\biggr)^2\\
 ={}&\frac{\abs{u}}2+\frac{2J}{9}\mu_1-\alpha^2\mu_1^2
 -\frac{2\abs{u}}9
 \biggl(\mu_2-\mu_3
       -\frac{3\mu_1\abs{\beta}}{\sqrt{2\abs{u}}}
       +\frac12\biggr)^2.
\end{aligned}
\end{equation*}
Conditions~\eqref{eq:quartic-Jpositive} and~\eqref{eq:quartic-weight3} give $\alpha\ne0$.
Maximize the quadratic in $\mu_1$ and set the last square to zero.
Together with $\mu_1+\mu_2+\mu_3=1$, this gives the weights $\mu^*=(\mu_1^*, \mu_2^*, \mu_3^*)$ defined by
\begin{equation}\label{eq:quartic-optimal-weights}
\begin{aligned}
 \mu_1^* \=\frac{J}{9\alpha^2}, \quad
 \mu_2^* \=\frac{9\alpha^2-2J(1-3\abs{\beta}/\sqrt{2\abs{u}})}{36\alpha^2},\quad
 \mu_3^* \=\frac{27\alpha^2-2J(1+3\abs{\beta}/\sqrt{2\abs{u}})}{36\alpha^2}.
\end{aligned}
\end{equation}
Conditions~\eqref{eq:quartic-Jpositive}--\eqref{eq:quartic-weight3} give $\mu_1^*,\,\mu_2^*,\,\mu_3^*>0$.
For this choice, $\Delta(\mu^*)=\frac{\abs{u}}{2}+\frac{J^2}{81\alpha^2}$, so $\Delta(\mu^*)-\frac34=\mathcal M\big/\bigl(162\abs{u}\alpha^2\bigr)\ge0$ by~\eqref{eq:quartic-Mnonnegative}.
Hence $R^2\ge Q_{\mu^*}(\abs{\gamma},\abs{\delta})\ge\Delta(\mu^*)\ge \frac34$.
Since $R^2\le\frac34$, we conclude that $R^2=\frac34$.
\end{proof}

\smallskip

\subsection{Proof of Theorem~\ref{prop:fourpoints}}
\begin{proof}
The radius bound follows from Lemma~\ref{lem:quartic-elimination}.

Now we suppose that $R^2=3/4$.
By the proof of Lemma~\ref{lem:quartic-elimination},
\begin{equation*}
 \frac34=R^2\ge Q_{\mu^*}(\abs{\gamma},\abs{\delta})
 \ge\Delta(\mu^*)=\frac34+\frac{\mathcal M}{162\abs{u}\alpha^2}\ge\frac34.
\end{equation*}
Thus $\mathcal M=0$, so $\lambda=1$ by Lemma~\ref{lem:quartic-algebra}.
By~\eqref{eq:quartic-center}, $u=-1$ and $p_t=4/3$, so we may take $\eta=0$.
Then $p_\tau=4/3$ by~\eqref{eq:quartic-rotation}.
Substituting these values into the definition of $J$ and~\eqref{eq:quartic-optimal-weights} gives $J=6$ and $(\mu_1^*, \mu_2^*, \mu_3^*)=(3/8, 1/16, 9/16)$.
By~\eqref{eq:quartic-quadratic}, $Q_{\mu^*}(\abs{\gamma},\abs{\delta})=\Delta(\mu^*)$ requires $\abs{\gamma}=1/2$, $\abs{\delta}=\sqrt2/6$, and $\abs{c_\tau}^2=11/36$.
Consequently, the second term in~\eqref{eq:quartic-radius} yields
\begin{equation*}
 \frac34\ge a_0^2+\frac29+\abs{c_\tau}^2+K
 \ge a_0^2+\frac29+\abs{c_\tau}^2+\frac{2\sqrt2}{3}\abs{\delta}
 =a_0^2+\frac34,
\end{equation*}
so $a_0=0$.
Here $p_\tau-c_\tau=c_w$ by~\eqref{eq:quartic-center} and~\eqref{eq:quartic-rotation}.
As $0\in\overline{D}(c_w,R)$, $\abs{p_\tau-c_\tau}\le R=\sqrt3/2$.
This excludes $\gamma=-1/2$, since then $\abs{p_\tau-c_\tau}^2=\bigl(-\frac12-\frac43\bigr)^2+\frac{1}{18}=\frac{41}{12}>\frac34=R^2$.
After complex conjugation if necessary, take $c_\tau\=\frac12+\frac{\sqrt2}{6}\I$.
Since $a_0=0$, $\tau_j$ are purely imaginary by~\eqref{eq:quartic-ellipse} for $j=2,\,3,\,4$. As $\eta=0$, it follows from~\eqref{eq:quartic-mean} that
$\frac13\sum_{j=2}^4\abs{\tau_j-c_\tau}^2
 =\frac49+\abs{c_\tau}^2=\frac34$,
which implies that $\abs{\tau_j-c_\tau}^2=\frac34$ for $j=2,\,3,\,4$. Since $\eta=0$, $\tau_j=t_j=\frac43-w_j$ by~\eqref{eq:quartic-center} and~\eqref{eq:quartic-rotation}. Since $\abs{\tau_j-c_\tau}^2=\frac34$ and $\sum_{j=2}^4\tau_j=0$, we obtain
\begin{equation*}
 \prod_{j=2}^4(w-w_j)
 =\biggl(w-\frac43+\frac{2\sqrt2}3\I\biggr)
  \biggl(w-\frac43-\frac{\sqrt2}3\I\biggr)^2.
\end{equation*}
Hence $w_2w_3w_4=(88+4\sqrt2\,\I)/27$. The value $(88-4\sqrt2\,\I)/27$ is obtained by conjugation. Conversely, for $\lambda=1$ and $w_2w_3w_4=(88+4\sqrt2\,\I)/27$,~\eqref{eq:quartic-sums} determines this polynomial. Its zeros and $0$ lie on $\Absbig{w-\frac56+\frac{\sqrt2}{6}\I}=\frac{\sqrt3}{2}$, so equality follows from the radius bound. The case $w_2w_3w_4=(88-4\sqrt2\,\I)/27$ follows by conjugation.
\end{proof}
\smallskip

\FloatBarrier
\appendix
\section*{Appendix}
\hypersetup{bookmarksdepth=2}
\renewcommand{\thesubsection}{\Alph{subsection}}
\renewcommand{\theHsubsection}{appendix.\Alph{subsection}}
\numberwithin{equation}{subsection} \numberwithin{table}{subsection} \numberwithin{theorem}{subsection}
\renewcommand{\theHtheorem}{\theHsubsection.\arabic{theorem}}
\makeatletter
\renewcommand{\subsection}{\@startsection{subsection}{2}%
  \z@{.7\linespacing\@plus.7\linespacing}{.3\linespacing}%
  {\normalfont\bfseries}}
\renewcommand{\subsubsection}{\@startsection{subsubsection}{3}%
  \z@{.5\linespacing\@plus.7\linespacing}{-.5em}%
  {\normalfont\itshape}}
\makeatother
\subsection{Positive polynomial expansions}\label{app:finiteidentities}
We verify the positivity of $P$ and $Q$ used in Lemma~\ref{lem:finiteexp} for all rows of Table~\ref{tab:finiteanalytic} here. Recall from the proof of Lemma~\ref{lem:finiteexp} that
\begin{equation*}
\begin{aligned}
P_{s,r_+,\lambda}(a)&=\bigl(1-\bigl(1-a^2\bigr)k^2-\lambda a\bigr)\bigl(2s+r_+\bigl(1-a^2\bigr)\bigr)+2 \quad\text{ and}\\
Q_{s,r_+,\lambda}(a)&=\bigl(S_6(-\lambda a)-k^2\bigr)\bigl(2s+r_+\bigl(1-a^2\bigr)\bigr)+2,
\end{aligned}
\end{equation*}
where $k=(s+1)/(2s)$, $S_6(-\lambda a)=\sum_{j=0}^{5}(-\lambda a)^j/j!$, and $s,\,r_+,\,\lambda$ are as in Table~\ref{tab:finiteanalytic}.

For $n\in\N_0$, the Bernstein basis on $[0,1]$ is $\Bbern_{n,i}(z)\=\binom ni z^i(1-z)^{n-i}$, $i=0,\,1,\,\ldots,\,n$ \cite{Farouki2012}. We will use the following elementary fact.

\begin{lemma}\label{lem:bernsteinpositive}
Let $n\in\N_0$ and let $F(a)=\sum_{h=0}^n p_ha^h\in\R[a]$. Let $a_-,\,a_+\in\R$ satisfy $a_-<a_+$. Then
\begin{equation}\label{eq:finiteBernstein}
 F(a_-+(a_+-a_-)z)=\sum_{i=0}^n b_i\Bbern_{n,i}(z),
\end{equation}
where
\begin{equation}\label{eq:bernsteinrule}
 b_i=\sum_{j=0}^i\frac{\binom ij}{\binom nj}(a_+-a_-)^j
       \sum_{h=j}^n p_h\binom hj a_-^{h-j},\qquad i=0,\,1,\,\ldots,\,n.
\end{equation}
For every $a\in[a_-,a_+]$, $F(a)\ge\min_{0\le i\le n}b_i$. In particular, $F>0$ on $[a_-,a_+]$ if all $b_i>0$.
\end{lemma}
\begin{proof}
Substitute $a\=a_-+(a_+-a_-)z$ and use the binomial identity $z^j=\sum_{i=j}^n\binom ij\Bbern_{n,i}(z) \big/\binom nj$. This gives \eqref{eq:finiteBernstein} with the coefficients in \eqref{eq:bernsteinrule}. The lower bound follows because $\Bbern_{n,i}(z)\ge0$ and $\sum_{i=0}^n\Bbern_{n,i}(z)=1$ for $z\in[0,1]$.
\end{proof}
\smallskip

For rows 1--5, the data below give the coefficients $b_i$ on each indicated interval.
We use $n=4$ for $P$ and $n=7$ for $Q$. For row~6, we record only the minimum coefficients for each $s$. All coefficients are computed by \eqref{eq:bernsteinrule} (see Appendix~\ref{app:code} for code). The intervals cover $[(s-1)/(s+1),1]$, so Lemma~\ref{lem:bernsteinpositive} gives the required positivity of both polynomials.

For rows 1--5, $s=2$ and the values of $r_+$ and $\lambda$ are given in Table~\ref{tab:finiteanalytic}. Tables~\ref{tab:BernsteinP} and~\ref{tab:BernsteinQ} list the coefficients for $P$ and $Q$, respectively.

\begin{appendixtable}
\begin{longtable}{@{}>{\centering\arraybackslash}m{28pt}>{\centering\arraybackslash}m{70pt}>{\centering\arraybackslash}m{\dimexpr\linewidth-126pt\relax}@{}}
\toprule
\textbf{Row}&{\boldmath$[a_-,a_+]$}&{\boldmath$(b_0,\,b_1,\,\ldots,\,b_n)$}\\
\midrule
\endfirsthead
\multicolumn{3}{@{}l@{}}{\footnotesize Table~\thetable\ (continued)}\\
\toprule
\textbf{Row}&{\boldmath$[a_-,a_+]$}&{\boldmath$(b_0,\,b_1,\,\ldots,\,b_n)$}\\
\midrule
\endhead
\midrule
\multicolumn{3}{@{}r@{}}{\footnotesize Continued on next page}\\
\endfoot
\bottomrule
\caption[Coefficients for $P$]{Coefficients for $P$: the coefficients are taken on $[a_-,a_+]$.}\label{tab:BernsteinP}\\
\endlastfoot
1&$[\frac{1}{3},\frac{2}{3}]$&$\dfrac{1}{1080}(3088,\,1975,\,1117,\,569,\,311)$\\*
\addlinespace[4pt]
1&$[\frac{2}{3},1]$&$\dfrac{1}{1080}(311,\,53,\,85,\,387,\,864)$\\
\addlinespace[6pt]
2&$[\frac{1}{3},\frac{2}{3}]$&$\dfrac{1}{72}(232,\,151,\,89,\,50,\,32)$\\
\addlinespace[4pt]
2&$[\frac{2}{3},1]$&$\dfrac{1}{72}(32,\,14,\,17,\,39,\,72)$\\
\addlinespace[6pt]
3&$[\frac{1}{3},1]$&$\dfrac{1}{180}(968,\,353,\,40,\,141,\,288)$\\
\addlinespace[6pt]
4&$[\frac{1}{3},1]$&$\dfrac{1}{108}(752,\,359,\,158,\,207,\,216)$\\*
\addlinespace[6pt]
5&$[\frac{1}{3},1]$&$\dfrac{1}{540}(20728,\,14353,\,10750,\,9621,\,1728)$\\
\addlinespace[4pt]
\end{longtable}
\end{appendixtable}

\begin{appendixtable}
\begin{longtable}{@{}>{\centering\arraybackslash}m{28pt}>{\centering\arraybackslash}m{70pt}>{\centering\arraybackslash}m{\dimexpr\linewidth-126pt\relax}@{}}
\toprule
\textbf{Row}&{\boldmath$[a_-,a_+]$}&{\boldmath$(b_0,\,b_1,\,\ldots,\,b_n)$}\\
\midrule
\endfirsthead
\multicolumn{3}{@{}l@{}}{\footnotesize Table~\thetable\ (continued)}\\
\toprule
\textbf{Row}&{\boldmath$[a_-,a_+]$}&{\boldmath$(b_0,\,b_1,\,\ldots,\,b_n)$}\\
\midrule
\endhead
\midrule
\multicolumn{3}{@{}r@{}}{\footnotesize Continued on next page}\\
\endfoot
\bottomrule
\caption[Coefficients for $Q$]{Coefficients for $Q$: the coefficients are taken on $[a_-,a_+]$.}\label{tab:BernsteinQ}\\
\endlastfoot
1&$[\frac{1}{3},1]$&$\dfrac{1}{45927000000}\begin{aligned}
 &\bigl(142664261971,\,92636217231,\,54085682241,\,27661100589,\,\\
 &\phantom{\bigl(}13237090551,\,10366710723,\,18464109813,\,37549256913\bigr)
\end{aligned}$\\*
\addlinespace[6pt]
2&$[\frac{1}{3},1]$&$\dfrac{1}{31352832}\begin{aligned}
 &\bigl(107187871,\,68782917,\,38986295,\,18482157,\,\\
 &\phantom{\bigl(}7325775,\,5291541,\,12020967,\,27530685\bigr)
\end{aligned}$\\
\addlinespace[6pt]
3&$[\frac{1}{3},\frac{2}{3}]$&$\dfrac{1}{15309000000}\begin{aligned}
 &\bigl(81244064167,\,65641437326,\,51632630429,\,39341716568,\,\\
 &\phantom{\bigl(}28857196408,\,20240312944,\,13531228672,\,8755319888\bigr)
\end{aligned}$\\
\addlinespace[4pt]
3&$[\frac{2}{3},1]$&$\dfrac{1}{15309000000}\begin{aligned}
 &\bigl(8755319888,\,3979411104,\,1136677808,\,252496296,\,\\
 &\phantom{\bigl(}1337456664,\,4393657107,\,9421115346,\,16426551897\bigr)
\end{aligned}$\\
\addlinespace[6pt]
4&$[\frac{1}{3},\frac{2}{3}]$&$\dfrac{1}{3674160}\begin{aligned}
 &\bigl(24200428,\,19946174,\,16081651,\,12646805,\,\\
 &\phantom{\bigl(}9672214,\,7181036,\,5190493,\,3713297\bigr)
\end{aligned}$\\
\addlinespace[4pt]
4&$[\frac{2}{3},1]$&$\dfrac{1}{3674160}\begin{aligned}
 &\bigl(3713297,\,2236101,\,1272252,\,834462,\,\\
 &\phantom{\bigl(}930933,\,1566783,\,2745414,\,4470228\bigr)
\end{aligned}$\\
\addlinespace[6pt]
5&$[\frac{1}{3},\frac{2}{3}]$&$\dfrac{1}{26244000000}\begin{aligned}
 &\bigl(873896393204,\,750898524862,\,634475839223,\,526215006496,\,\\
 &\phantom{\bigl(}427486552136,\,339474732928,\,263202251264,\,199552311856\bigr)
\end{aligned}$\\*
\addlinespace[4pt]
5&$[\frac{2}{3},1]$&$\dfrac{1}{26244000000}\begin{aligned}
 &\bigl(199552311856,\,135902372448,\,84874975296,\,47353325112,\,\\
 &\phantom{\bigl(}24081128928,\,15684652809,\,22692121002,\,45552961764\bigr)
\end{aligned}$\\
\addlinespace[4pt]
\end{longtable}
\end{appendixtable}

For row 6, $3\le s\le8$, $\lambda=13/12$. For each such $s$, we take $r_+\=r_0(s)-1$ as specified in Table~\ref{tab:finiteanalytic} and divide $[(s-1)/(s+1),1]$ into $m$ equal subintervals, with $m$ chosen separately for $P$ and $Q$ as listed in Table~\ref{tab:mergedBernstein}. Explicitly, their endpoints are
\begin{equation*}
 a_-\=\frac{s-1}{s+1}+\frac{2j}{m(s+1)}\quad\text{ and }\quad
 a_+\=a_-+\frac{2}{m(s+1)},\qquad j=0,\,1,\,\ldots,\,m-1.
\end{equation*}
Table~\ref{tab:mergedBernstein} gives the exact minimum of the coefficients in \eqref{eq:bernsteinrule} over all indices and subintervals.

\begin{appendixtable}
\begin{longtable}{@{}>{\centering\arraybackslash}m{28pt}>{\centering\arraybackslash}m{24pt}>{\centering\arraybackslash}m{\dimexpr(\linewidth-132pt)/2\relax}>{\centering\arraybackslash}m{24pt}>{\centering\arraybackslash}m{\dimexpr(\linewidth-132pt)/2\relax}@{}}
\toprule
\multirow{2}{*}{\boldmath$s$}&\multicolumn{2}{c}{\boldmath$P$}&\multicolumn{2}{c}{\boldmath$Q$}\\
\cmidrule(lr){2-3}\cmidrule(lr){4-5}
&{\boldmath$m$}&{\boldmath$\min b_i$}&{\boldmath$m$}&{\boldmath$\min b_i$}\\
\midrule
\endfirsthead
\multicolumn{5}{@{}l@{}}{\footnotesize Table~\thetable\ (continued)}\\
\toprule
\multirow{2}{*}{\boldmath$s$}&\multicolumn{2}{c}{\boldmath$P$}&\multicolumn{2}{c}{\boldmath$Q$}\\
\cmidrule(lr){2-3}\cmidrule(lr){4-5}
&{\boldmath$m$}&{\boldmath$\min b_i$}&{\boldmath$m$}&{\boldmath$\min b_i$}\\
\midrule
\endhead
\midrule
\multicolumn{5}{@{}r@{}}{\footnotesize Continued on next page}\\
\endfoot
\bottomrule
\caption[Bernstein coefficients]{Bernstein coefficients: $m$ is the number of equal subintervals in $a$.}\label{tab:mergedBernstein}\\
\endlastfoot
3&3&$\dfrac{2483}{31104}$&5&$\dfrac{3144211937}{460800000000}$\\*
\addlinespace[6pt]
4&2&$\dfrac{223}{800}$&1&$\dfrac{30949417}{29030400}$\\
\addlinespace[6pt]
5&1&$\dfrac{29}{60}$&1&$\dfrac{26355163}{14929920}$\\
\addlinespace[6pt]
6&1&$\dfrac{4}{21}$&1&$\dfrac{4864607}{2488320}$\\
\addlinespace[6pt]
7&2&$\dfrac{5293}{25088}$&1&$\dfrac{223643183}{104509440}$\\*
\addlinespace[6pt]
8&3&$\dfrac{13255}{472392}$&1&$\dfrac{4333247}{1866240}$\\
\addlinespace[4pt]
\end{longtable}
\end{appendixtable}

\subsection{Rational comparisons for the degree ranges}\label{app:finitearithmetic}
We verify the bounds in Table~\ref{tab:finiteanalytic} by direct substitution in the following three cases. For convenience, we recall the corresponding upper bounds before substituting the data. The code is given in Appendix~\ref{app:code}.

\subsubsection{Bounds for $E$}
By \eqref{eq:finitee},
\begin{equation*}
 E\le\frac{4k^2}{\lambda^2}
       \sum_{j=0}^{s-1}(j+1)
       \biggl(\frac{2k}{\lambda}\biggr)^j
       \frac{r_-^{j+1}}{(r_--1)^{j+2}}.
\end{equation*}
Substituting $k=(s+1)/(2s)$ and the tabulated values of $s$, $\lambda$, and $r_-$ gives Table~\ref{tab:integralarithmetic}.
\begin{appendixtable}
\begin{longtable}{@{}>{\centering\arraybackslash}m{28pt}>{\centering\arraybackslash}m{28pt}>{\centering\arraybackslash}m{\dimexpr\linewidth-84pt\relax}@{}}
\toprule
\textbf{Row}&{\boldmath$s$}&\textbf{Numerical comparison}\\
\midrule
\endfirsthead
\multicolumn{3}{@{}l@{}}{\footnotesize Table~\thetable\ (continued)}\\
\toprule
\textbf{Row}&{\boldmath$s$}&\textbf{Numerical comparison}\\
\midrule
\endhead
\midrule
\multicolumn{3}{@{}r@{}}{\footnotesize Continued on next page}\\
\endfoot
\bottomrule
\caption[Integral comparisons]{Integral comparisons.}\label{tab:integralarithmetic}\\
\endlastfoot
1&2&$\displaystyle\frac{250}{1521}\biggl(1+\frac{100}{39}\biggr)<\frac{59}{100}$\\*
\addlinespace[6pt]
2&2&$\displaystyle\frac{99}{625}\biggl(1+\frac{66}{25}\biggr)<\frac{58}{100}$\\
\addlinespace[6pt]
3&2&$\displaystyle\frac{325}{1936}\biggl(1+\frac{65}{22}\biggr)<\frac{67}{100}$\\
\addlinespace[6pt]
4&2&$\displaystyle\frac{25}{256}\biggl(1+\frac{25}{8}\biggr)<\frac{41}{100}$\\
\addlinespace[6pt]
5&2&$\displaystyle\frac{6750}{41209}\biggl(1+\frac{900}{203}\biggr)<\frac{90}{100}$\\
\addlinespace[6pt]
6&3&$\displaystyle\frac{17}{169}\sum\nolimits_{j=0}^{2}(j+1)\biggl(\frac{17}{13}\biggr)^j<\frac{88}{100}$\\
\addlinespace[4pt]
6&4&$\displaystyle\frac{18}{325}\sum\nolimits_{j=0}^{3}(j+1)\biggl(\frac65\biggr)^j<\frac{88}{100}$\\
\addlinespace[4pt]
6&5&$\displaystyle\frac{148}{4225}\sum\nolimits_{j=0}^{4}(j+1)\biggl(\frac{74}{65}\biggr)^j<\frac{88}{100}$\\
\addlinespace[4pt]
6&6&$\displaystyle\frac{200}{8281}\sum\nolimits_{j=0}^{5}(j+1)\biggl(\frac{100}{91}\biggr)^j<\frac{88}{100}$\\
\addlinespace[4pt]
6&7&$\displaystyle\frac{45}{2548}\sum\nolimits_{j=0}^{6}(j+1)\biggl(\frac{15}{14}\biggr)^j<\frac{88}{100}$\\*
\addlinespace[4pt]
6&8&$\displaystyle\frac{41}{3042}\sum\nolimits_{j=0}^{7}(j+1)\biggl(\frac{41}{39}\biggr)^j<\frac{88}{100}$\\
\addlinespace[4pt]
\end{longtable}
\end{appendixtable}

\subsubsection{Bounds for $M$ in rows 4--6}
By \eqref{eq:finitemendpoint},
\begin{equation*}
 M\le\sum_{j=0}^{s-1}\frac{k^j}{j!}\,r_-^jk^{r_-}
       +\frac{k^s}{(s-1)!}\frac{r_-^sk^{r_-}}{r_-+s}.
\end{equation*}
Substituting $k=(s+1)/(2s)$ and the tabulated values of $s$ and $r_-$ gives Table~\ref{tab:endpointarithmetic}.
\begin{appendixtable}
\begin{longtable}{@{}>{\centering\arraybackslash}m{28pt}>{\centering\arraybackslash}m{28pt}>{\centering\arraybackslash}m{\dimexpr\linewidth-84pt\relax}@{}}
\toprule
\textbf{Row}&{\boldmath$s$}&\textbf{Numerical comparison}\\
\midrule
\endfirsthead
\multicolumn{3}{@{}l@{}}{\footnotesize Table~\thetable\ (continued)}\\
\toprule
\textbf{Row}&{\boldmath$s$}&\textbf{Numerical comparison}\\
\midrule
\endhead
\midrule
\multicolumn{3}{@{}r@{}}{\footnotesize Continued on next page}\\
\endfoot
\bottomrule
\caption[Endpoint comparisons]{Endpoint comparisons.}\label{tab:endpointarithmetic}\\
\endlastfoot
4&2&$\displaystyle\biggl(1+\frac{75}{4}+\frac{(75/4)^2}{27}\biggr)\biggl(\frac34\biggr)^{25}<\frac{3}{100}$\\*
\addlinespace[6pt]
5&2&$\displaystyle\biggl(1+\frac{45}{2}+\frac{(45/2)^2}{32}\biggr)\biggl(\frac34\biggr)^{30}<\frac{1}{100}$\\
\addlinespace[6pt]
6&3&$\displaystyle\biggl(\sum\nolimits_{j=0}^{2}\frac{(34/3)^j}{j!}+\frac{(34/3)^3}{2!\cdot20}\biggr)\biggl(\frac23\biggr)^{17}<\frac{12}{100}$\\
\addlinespace[4pt]
6&4&$\displaystyle\biggl(\sum\nolimits_{j=0}^{3}\frac{(65/4)^j}{j!}+\frac{(65/4)^4}{3!\cdot30}\biggr)\biggl(\frac58\biggr)^{26}<\frac{12}{100}$\\
\addlinespace[4pt]
6&5&$\displaystyle\biggl(\sum\nolimits_{j=0}^{4}\frac{(111/5)^j}{j!}+\frac{(111/5)^5}{4!\cdot42}\biggr)\biggl(\frac35\biggr)^{37}<\frac{12}{100}$\\
\addlinespace[4pt]
6&6&$\displaystyle\biggl(\sum\nolimits_{j=0}^{5}\frac{(175/6)^j}{j!}+\frac{(175/6)^6}{5!\cdot56}\biggr)\biggl(\frac7{12}\biggr)^{50}<\frac{12}{100}$\\
\addlinespace[4pt]
6&7&$\displaystyle\biggl(\sum\nolimits_{j=0}^{6}\frac{(260/7)^j}{j!}+\frac{(260/7)^7}{6!\cdot72}\biggr)\biggl(\frac47\biggr)^{65}<\frac{12}{100}$\\*
\addlinespace[4pt]
6&8&$\displaystyle\biggl(\sum\nolimits_{j=0}^{7}\frac{(369/8)^j}{j!}+\frac{(369/8)^8}{7!\cdot90}\biggr)\biggl(\frac9{16}\biggr)^{82}<\frac{12}{100}$\\
\addlinespace[4pt]
\end{longtable}
\end{appendixtable}

\subsubsection{Bounds for $M$ in rows 1--3}
By \eqref{eq:finitemsharp},
\begin{equation*}
 M\le\biggl(1+\frac{3(r_-+1)}4\biggr)
       \biggl(\frac9{16}-\frac{2}{4+8r_+/9}\biggr)^{\lfloor r_-/2\rfloor}
       +\frac{(3r_-/4)^2(3/4)^{r_-}}{r_-+2}.
\end{equation*}
Substituting $(r_-,r_+)=(10,10)$, $(11,12)$, and $(13,24)$ gives Table~\ref{tab:sharpendpointarithmetic}.
\begin{appendixtable}
\begin{longtable}{@{}>{\centering\arraybackslash}m{28pt}>{\centering\arraybackslash}m{28pt}>{\centering\arraybackslash}m{\dimexpr\linewidth-84pt\relax}@{}}
\toprule
\textbf{Row}&{\boldmath$s$}&\textbf{Numerical comparison}\\
\midrule
\endfirsthead
\multicolumn{3}{@{}l@{}}{\footnotesize Table~\thetable\ (continued)}\\
\toprule
\textbf{Row}&{\boldmath$s$}&\textbf{Numerical comparison}\\
\midrule
\endhead
\midrule
\multicolumn{3}{@{}r@{}}{\footnotesize Continued on next page}\\
\endfoot
\bottomrule
\caption[Sharper endpoint comparisons]{Sharper endpoint comparisons.}\label{tab:sharpendpointarithmetic}\\
\endlastfoot
1&2&$\displaystyle\frac{37}{4}\biggl(\frac{189}{464}\biggr)^5+\frac{75}{16}\biggl(\frac34\biggr)^{10}<\frac{37}{100}$\\*
\addlinespace[6pt]
2&2&$\displaystyle10\biggl(\frac{75}{176}\biggr)^5+\frac{1089}{208}\biggl(\frac34\biggr)^{11}<\frac{37}{100}$\\*
\addlinespace[6pt]
3&2&$\displaystyle\frac{23}{2}\biggl(\frac{147}{304}\biggr)^6+\frac{507}{80}\biggl(\frac34\biggr)^{13}<\frac{30}{100}$\\
\addlinespace[4pt]
\end{longtable}
\end{appendixtable}

\subsection{Polynomial estimates for Lemma~\ref{lem:quartic-algebra}}\label{app:quartic}
\begingroup
\setlength{\parskip}{2pt plus 1pt}
\postdisplaypenalty=10000
\clubpenalty=10000
\widowpenalty=10000
We use the notation of Lemma~\ref{lem:quartic-algebra} and its proof, where $\kappa=x+\I y$, $\abs{\kappa}\ge1$, $1\le\abs{\lambda}<3/2$, and $0<\abs{u}<27/16$.

By~\eqref{eq:quartic-joukowski}, if $x\ge2/3$, then $64\abs{\kappa}^2\abs{\lambda}^2=\bigl(9\abs{\kappa}^2+6x+1\bigr)\bigl(\abs{\kappa}^2+6x+9\bigr)\ge\bigl(9\abs{\kappa}^2+5\bigr)\bigl(\abs{\kappa}^2+13\bigr)>144\abs{\kappa}^2$. This contradicts $\abs{\lambda}<3/2$, so $x<2/3$.

Set $t\=\abs{\kappa}^2-1$ and $v\=\abs{\kappa}^2+6x+1$,
and
\begin{align*}
 D_1&\=36y^2=-t^2+2tv+32t-v^2+4v+32,\\
 D_2&\=64\abs{\kappa}^2\bigl(\tfrac94-\abs{\lambda}^2\bigr)
       =-8tv+80t-v^2-16v+80,\\
 D_3&\=32\abs{\kappa}^2\bigl(\tfrac{27}{16}-\abs{u}\bigr)
       =-8t^2-2tv+22t+v^2-4v+22.
\end{align*}
Then $t,\,v\ge0$, $6+t-v=4-6x>0$, $D_1\ge0$, and $D_2,\,D_3>0$.

Substituting $\abs{\kappa}^2=1+t$ and $x=(v-t-2)/6$ gives the identities below (see Appendix~\ref{app:code} for code), in which each $\sigma_j$ has nonnegative coefficients and a positive constant term.

\subsubsection{Positivity of $\mathcal N_1$, $\mathcal N_2$, and $\mathcal N_3$}
\begin{equation}\label{eq:quartic-C1}
 18\mathcal N_1=6D_1+6D_3+\sigma_1,
\end{equation}
where $\sigma_1\=6t^{2}+6tv+60t+6v^{2} +48v+60$.
\begin{equation}\label{eq:quartic-C2}
\begin{aligned}
 81\mathcal N_2={}&\bigl(\tfrac{21}{16}tv^2+\tfrac12v^3\bigr)D_1
 +\bigl(\tfrac1{16}t^3+\tfrac{55}8tv\bigr)D_3\\
 &+\bigl(\tfrac{25}{16}t^3+\tfrac5{16}tv^2+\tfrac{115}8tv
 +\tfrac{11}8v^2+\tfrac{23}2v+\tfrac{409}{16}\bigr)D_2+\sigma_2,
\end{aligned}
\end{equation}
where $\sigma_2\=\tfrac{1}{2}t^{5}+\tfrac{101}{8}t^{4}v+\tfrac{13}{8}t^{4}+\tfrac{45}{16}t^{3}v^{2} +\tfrac{1}{4}t^{3}v+\tfrac{7181}{8}t^{3}+\tfrac{3}{8}t^{2}v^{3}+\tfrac{7}{4}t^{2}v^{2} +\tfrac{4395}{4}t^{2}v+3072t^{2}+\tfrac{5}{8}tv^{4}+\tfrac{49}{4}tv^{3} +\tfrac{657}{2}tv^{2}+\tfrac{3069}{4}tv+2051t+\tfrac{1}{2}v^{5} +\tfrac{11}{8}v^{4}+\tfrac{3}{2}v^{3}+\tfrac{57}{16}v^{2}+v +3$.
\begin{equation}\label{eq:quartic-C3}
\begin{aligned}
 \mathcal N_3= \bigl(\tfrac{23}8tv+\tfrac{41}8v+\tfrac{89}8\bigr)(6+t-v)
  +\tfrac16v^2D_1+\tfrac23vD_3 
  +\bigl(\tfrac7{12}tv+\tfrac72t+\tfrac18v^2+\tfrac23v+2\bigr)D_2
  +\sigma_3,
\end{aligned}
\end{equation}
where $\sigma_3\=\tfrac{64}{9}t^{3}v+\tfrac{512}{9}t^{3}+\tfrac{47}{18}t^{2}v^{2}+\tfrac{113}{72}t^{2}v +\tfrac{40}{9}t^{2}+\tfrac{13}{36}tv^{3}+\tfrac{91}{72}tv^{2}+\tfrac{181}{72}tv +\tfrac{287}{72}t+\tfrac{7}{24}v^{4}+\tfrac{4}{9}v^{3}+\tfrac{113}{72}v^{2} +\tfrac{91}{72}v+\tfrac{29}{36}$. 

Hence, $\mathcal N_1,\,\mathcal N_2,\,\mathcal N_3>0$.

\subsubsection{Bounds for $H_1$ and $W$}
\begin{equation}\label{eq:quartic-C6}
 H_1-12\abs{\kappa}^2=\tfrac12v(6+t-v)+\tfrac12D_3+\sigma_4,
\end{equation}
where $\sigma_4\=\tfrac{28}{9}t^{2}+\tfrac{41}{18}tv+\tfrac{49}{9}t+\tfrac{1}{9}v^{2} +\tfrac{23}{9}v+\tfrac{49}{9}$. Thus $H_1-12\abs{\kappa}^2>0$.
Next,
\begin{equation}\label{eq:quartic-C4}
 36W=41D_1+\sigma_5,
\end{equation}
where $\sigma_5\=201t^{2}+30tv+48t+21v^{2} +60v+48$. Hence $W>0$.
Finally,
\begin{equation}\label{eq:quartic-C5}
\begin{aligned}
 W^2-48H_1\abs{\kappa}^2= \tfrac{13}3v(6+t-v)+\tfrac13v^2D_1 
  +\bigl(2tv+\tfrac23t+\tfrac13v^2+\tfrac{10}3v+\tfrac23\bigr)D_2
   +\sigma_6,
\end{aligned}
\end{equation}
where $\sigma_6\=\tfrac{1600}{81}t^{4}+\tfrac{2240}{81}t^{3}v+\tfrac{30656}{81}t^{3}+\tfrac{569}{27}t^{2}v^{2} +\tfrac{1360}{27}t^{2}v+\tfrac{10448}{27}t^{2}+\tfrac{44}{81}tv^{3}+\tfrac{362}{27}tv^{2} +\tfrac{923}{27}tv+\tfrac{1376}{81}t+\tfrac{79}{81}v^{4}+\tfrac{34}{81}v^{3} +\tfrac{335}{27}v^{2}+\tfrac{1414}{81}v+\tfrac{688}{81}$. Hence $W>4\sqrt{3H_1\abs{\kappa}^2}$.

\subsubsection{A lower bound for $\mathcal M$}
Substituting~\eqref{eq:quartic-joukowski-invariants} and $J=9\mathcal N_1\big/\bigl(32\abs{\kappa}^2\bigr)$ into the definition of $\mathcal M$ gives
\begin{equation}\label{eq:quartic-master}
 \mathcal M=\tfrac{81}{4096}\abs{\kappa}^{-4}
 \bigl(256\bigl(\abs{\kappa}^4-1\bigr)^2+v\mathcal G\bigr),
\end{equation}
where
$ \mathcal G\= \frac13\bigl( (6+t-v)\bigl(40t^3+24t^2v+80tv+288t+8v^2+192v+512\bigr)   
 +24t^4+464t^3+4t^2v+2400t^2+2tv^3
 +2880t+v^4+1024\bigr)$.
Thus $\mathcal G\ge1024/3$, and
\begin{equation}\label{eq:quartic-Mcertificate}
 \mathcal M-
 \tfrac{1}{16}\abs{\kappa}^{-4} \bigl( 81(\abs{\kappa}^4-1)^2+108v \bigr)
 =\tfrac{81}{4096}v\abs{\kappa}^{-4}
   \bigl(\mathcal G-\tfrac{1024}{3}\bigr)\ge0.
\end{equation}
Thus $\mathcal M\ge0$, and equality implies $\abs{\kappa}^2=1$ and $v=0$, hence $x=-1/3$.

\par\endgroup
\subsection{Code for Appendices A--C}\label{app:code}
The following two programs use exact rational arithmetic and run with Python~3.9+.
\subsubsection{Tables in Appendices A and B}
The parameters of Table~\ref{tab:finiteanalytic} and the comparison bounds are inputs.
The code computes the coefficients and majorants in Appendices~\ref{app:finiteidentities} and~\ref{app:finitearithmetic}, with decimals only for displaying the majorants.
\begingroup
\lstset{language=Python,basicstyle=\linespread{1}\fontsize{8}{9}\selectfont\ttfamily,
 columns=fullflexible,keepspaces=true,showstringspaces=false,
 upquote=true,breaklines=false,aboveskip=6pt,belowskip=6pt}
\begin{lstlisting}[belowskip=0pt]
from fractions import Fraction as F
from itertools import zip_longest
from math import comb, factorial, lcm
CASES = [(2, 10, 10, F(13, 10)), (2, 11, 12, F(5, 4)),
         (2, 13, 24, F(11, 10)), (2, 25, 29, F(1)),
         (2, 30, 149, F(7, 10))]
CASES += [(s, (s + 1)**2 + 1, hi, F(13, 12))
          for s, hi in zip(range(3, 9), [119, 109, 71, 97, 127, 161])]
def add(p, q): return [a + b for a, b in zip_longest(p, q, fillvalue=0)]

def mul(p, q):
    out = [F(0)] * (len(p) + len(q) - 1)
    for i, a in enumerate(p):
        for j, b in enumerate(q):
            out[i + j] += a * b
    return out

def polynomials(s, r, lam):
    k, factor = F(s + 1, 2*s), [2*s + r, 0, -r]
    q = [(-lam)**j / factorial(j) for j in range(6)]
    q[0] -= k*k
    return (add(mul([1 - k*k, -lam, k*k], factor), [2]),
            add(mul(q, factor), [2]))

def bernstein(p, a, b):
    n = len(p) - 1
    c = [sum(p[h]*comb(h, j)*a**(h-j)*(b-a)**j
             for h in range(j, n+1)) for j in range(n+1)]
    return [sum(c[j]*F(comb(i, j), comb(n, j))
                for j in range(i+1)) for i in range(n+1)]

def minimum(p, s, m):
    a, d = F(s-1, s+1), F(2, m*(s+1))
    return min(min(bernstein(p, a+j*d, a+(j+1)*d)) for j in range(m))

def e_bound(s, r, lam):
    k = F(s+1, 2*s)
    return 4*r*k*k/(lam*lam*(r-1)**2) * sum(
        (j+1)*(2*r*k/(lam*(r-1)))**j for j in range(s))

def m_bound(s, lo, hi, sharp):
    k, y = F(s+1, 2*s), F(lo*(s+1), 2*s)
    if sharp:
        q = F(9, 16) - 2/(4 + F(8*hi, 9))
        return ((1 + F(3*(lo+1), 4))*q**(lo//2)
                + F(3*lo, 4)**2*F(3, 4)**lo/(lo+2))
    return k**lo * (sum(y**j/factorial(j) for j in range(s))
                    + y**s/(factorial(s-1)*(lo+s)))

print("\nTables A.1-A.2: row P/Q interval factor numerators")
for row, (s, lo, hi, lam) in enumerate(CASES[:5], 1):
    for name, p in zip("PQ", polynomials(s, hi, lam)):
        cuts = [F(1, 3), F(1)]
        if min(bernstein(p, *cuts)) <= 0:
            cuts.insert(1, F(2, 3))
        for a, b in zip(cuts, cuts[1:]):
            values = bernstein(p, a, b)
            d = lcm(*(v.denominator for v in values))
            print(row, name, f"[{a},{b}]", f"1/{d}",
                  tuple(int(d*v) for v in values))
\end{lstlisting}
\begin{lstlisting}[aboveskip=0pt]
print("\nTable A.3: s m_P min_P m_Q min_Q")
for s, lo, hi, lam in CASES[5:]:
    result = []
    for p in polynomials(s, hi, lam):
        m = 1
        while minimum(p, s, m) <= 0:
            m += 1
        result.extend([m, minimum(p, s, m)])
    print(s, *result)
# Chosen comparison bounds in Appendix B (hundredths).
bounds = [(59, 37), (58, 37), (67, 30), (41, 3), (90, 1), (88, 12)]
print("\nTables B.1-B.3: row s E-majorant < bound; M-majorant < bound")
for i, (s, lo, hi, lam) in enumerate(CASES, 1):
    row = min(i, 6)
    e, m = e_bound(s, lo, lam), m_bound(s, lo, hi, i <= 3)
    eb, mb = bounds[row-1]
    assert e < F(eb, 100) and m < F(mb, 100)
    print(row, s, f"{float(e):.12g} < {eb}/100;",
          f"{float(m):.12g} < {mb}/100")
\end{lstlisting}
\endgroup
\subsubsection{Identities in Appendix C}
The code computes $D_1,\,D_2,\,D_3,\,\sigma_1,\,\sigma_2,\,\ldots,\,\sigma_6$ and checks \eqref{eq:quartic-master} and~\eqref{eq:quartic-Mcertificate} after clearing denominators.

\begin{lstlisting}[language=Python,
 basicstyle=\linespread{1}\fontsize{8}{9}\selectfont\ttfamily,
 columns=fullflexible,keepspaces=true,showstringspaces=false,
 upquote=true,breaklines=false,aboveskip=6pt,belowskip=6pt]
from fractions import Fraction as F
class P(dict):
    def __init__(self, q=0):
        q = q if isinstance(q, dict) else {(0, 0): q}
        super().__init__((k, F(c)) for k, c in q.items() if c)
    def __add__(p, q):
        q = P(q)
        return P({k: p.get(k, 0)+q.get(k, 0) for k in p.keys()|q.keys()})
    __radd__ = __add__
    def __neg__(p): return P({k: -c for k, c in p.items()})
    def __sub__(p, q): return p + -P(q)
    def __rsub__(p, q): return P(q) + -p
    def __mul__(p, q):
        r = {}
        for (a, b), c in p.items():
            for (i, j), d in P(q).items():
                k = (a+i, b+j)
                r[k] = r.get(k, 0)+c*d
        return P(r)
    __rmul__ = __mul__
    def __truediv__(p, q): return p*(1/F(q))
    def __pow__(p, n):
        r = P(1)
        for _ in range(n): r = r*p
        return r

t, v = P({(1, 0): 1}), P({(0, 1): 1})
k2, x = 1+t, (v-t-2)/6  # k2 = |kappa|^2
h1 = (2*x*x-k2+6*x+1)**2+4*(x+3)**2*(k2-x*x)
h2 = (2*x*x-k2-1)**2+4*x*x*(k2-x*x)
a0, b02 = t*v, 4*(k2-x*x)*(3*k2+2*x+3)**2
n1 = h1-2*h2-a0
n2, n3 = b02-n1*(h2-a0), 3*b02-n1*(h2+a0)
d1 = 36*(k2-x*x)
d2 = 144*k2-(9*k2+6*x+1)*(k2+6*x+9)
d3, w = 9*(6*k2-h2), h1+6*h2-12*k2
assert h1*h2 == a0*a0+b02
sigma = [18*n1-6*d1-6*d3,
    81*n2-(F(21,16)*t*v*v+v**3/2)*d1-(t**3/16+F(55,8)*t*v)*d3
    -(F(25,16)*t**3+F(5,16)*t*v*v+F(115,8)*t*v+F(11,8)*v*v
      +F(23,2)*v+F(409,16))*d2,
    n3-(F(23,8)*t*v+F(41,8)*v+F(89,8))*(6+t-v)-v*v*d1/6
    -2*v*d3/3-(F(7,12)*t*v+F(7,2)*t+v*v/8+2*v/3+2)*d2,
    h1-12*k2-v*(6+t-v)/2-d3/2, 36*w-41*d1,
    w*w-48*h1*k2-F(13,3)*v*(6+t-v)-v*v*d1/3
    -(2*t*v+2*t/3+v*v/3+10*v/3+F(2,3))*d2]
for s in sigma:
    assert min(s.values()) >= 0 and s.get((0, 0), 0) > 0
l = 40*t**3+24*t*t*v+80*t*v+288*t+8*v*v+192*v+512
r = 24*t**4+464*t**3+4*t*t*v+2400*t*t+2*t*v**3+2880*t+v**4+1024
g = ((6+t-v)*l+r)/3
# phi = 4096*k2**3*M/729 after substituting the invariant formulas.
phi = h2*n1*n1/4-3*b02*(16*k2-3*h2)/4
# Clearing denominators in the two identities for M gives:
assert 9*phi == 256*k2*(k2*k2-1)**2+k2*v*g
assert 729*phi-256*k2*(81*(k2*k2-1)**2+108*v) == 81*k2*v*(g-F(1024,3))
assert min(l.values()) >= 0 and min((r-1024).values()) >= 0
for name, p in [("D1", d1), ("D2", d2), ("D3", d3)] + [
        (f"sigma{i}", s) for i, s in enumerate(sigma, 1)] + [("G", g)]:
    print(name, "=", " + ".join("*".join([str(c)]+[
          f"{z}^{n}" for z, n in zip("tv", powers) if n])
          for powers, c in sorted(p.items(), reverse=True)).replace("+ -", "- "))
\end{lstlisting}

\end{document}